\documentclass[reqno,11pt]{amsart}  

\usepackage{tikz}
\usetikzlibrary{decorations.pathreplacing,calligraphy}
\usepackage{todonotes}
\usepackage{gensymb}
\usepackage{amsmath}
\usepackage{amssymb} 
\usepackage{mathtools}   
\usepackage{graphicx}
\usepackage{graphics}  
\usepackage{color}
\usepackage{verbatim} 
\usepackage[normalem]{ulem} 
\usepackage[mathscr]{eucal}
\usepackage{upgreek}
\usepackage{enumerate} 
\usepackage{cite}
\usepackage{amsmath}
\usepackage{amsfonts}
\usepackage{amssymb}
\usepackage{bbm}
\usepackage{cancel}
\usepackage{graphicx}
\usepackage{multicol}
\usepackage[utf8]{inputenc}
\usepackage{framed}
\usepackage{setspace}
\usepackage{amsthm}
\usepackage{hyperref}
\usepackage{caption}
\usepackage{subcaption}
\usepackage{bbm}
\usepackage{graphicx}
\usepackage{svg}
\usepackage[titletoc]{appendix}
\numberwithin{equation}{section}  
 \usepackage{placeins}
 \usepackage{caption}
\usepackage[refpage,noprefix]{nomencl}
\usepackage{nomencl} 

\newtheorem{theorem}{Theorem}[section] 
\newtheorem{lemma}[theorem]{Lemma} 
\newtheorem{proposition}[theorem] {Proposition} 
\newtheorem{cor}[theorem]  {Corollary} 
\newtheorem{remark}[theorem]  {Remark} 
\newtheorem{definition}[theorem] {Definition}

\theoremstyle{definition}

 \usepackage{float}
\restylefloat{figure}

\DeclareMathAlphabet{\mathpzc}{OT1}{pzc}{m}{it}

\newcommand{\Babs}[1]{{\Bigl\lvert #1\Bigr\rvert}}

\newcommand{\abs}[1]{\left| #1 \right|}

\renewcommand{\L} {\Lambda} %

\def\a{\alpha}

\def\d{\delta} 
\newcommand{\e} {\varepsilon}

\def\l{\lambda}

\newfam\Bbbfam 
\font\tenBbb=msbm10 
\font\sevenBbb=msbm7 
\font\fiveBbb=msbm5 
\textfont\Bbbfam=\tenBbb 
\scriptfont\Bbbfam=\sevenBbb 
\scriptscriptfont\Bbbfam=\fiveBbb 
 
\newcommand{\B}     {\mathbb{B}} 
  
\newcommand{\R}     {\mathbb{R}} 
\newcommand{\Z}     {\mathbb{Z}} 
\newcommand{\N}     {\mathbb{N}} 
\renewcommand{\P}   {\mathbb{P}} 
 
\newcommand{\E}     {\mathbb{E}}

 \newcommand{\floor}[1]{\left\lfloor #1 \right\rfloor}

\def\1{{\mathchoice {1\mskip-4mu\mathrm l}      
{1\mskip-4mu\mathrm l} 
{1\mskip-4.5mu\mathrm l} {1\mskip-5mu\mathrm l}}} 
\newcommand{\ssup}[1] {{\scriptscriptstyle{({#1}})}} 
\def\comment#1{} 
\newtheoremstyle{thm}{2ex}{2ex}{\itshape\rmfamily}{} 
{\bfseries\rmfamily}{}{1.7ex}{} 
 
\newtheoremstyle{rem}{1.3ex}{1.3ex}{\rmfamily}{} 
{\itshape\rmfamily}{}{1.5ex}{} 
 
\newcommand{\bB} {\boldsymbol{B}}

\newcommand{\sfrak}{{\mathfrak{s}}}

\newcommand{\Ecal}   {{\mathcal E }}

\newcommand{\brho}{\boldsymbol{\rho}}

\newcommand{\Lcal}   {{\mathcal L }}

\newcommand{\Ocal}   {{\mathcal O }}

\newcommand{\Rcal}   {{\mathcal R }}

\newcommand{\Xcal}   {{\mathcal X }}

\newcommand{\pmf}{{\mathrm{MF}}}

 \newcommand{\ex}{{\rm e}} 
 
\renewcommand{\d}{{\rm d}}

\newcommand{\supp}{{\operatorname {supp}}}

\newcommand{\dist}{{\operatorname {dist}}}

\newcommand{\Exp}{\mathscr{E}\kern-0.2mm{\operatorname{xp}}}
\newcommand{\Log}{\mathscr{L}\kern-0.2mm{\operatorname{og}}}

\renewcommand{\emptyset} {\varnothing}

\newcommand{\Hpfm}{H^{\textsf{MF}}}

\newcommand\NoBlackBoxes{\global\overfullrule0pt}
\NoBlackBoxes

\newcommand{\rhoc}{{\rho_{\mathrm{c}}}}

\newcommand\mycom[2]{\genfrac{}{}{0pt}{}{#1}{#2}}

\newcommand{\Mfrak}{\mathrm{M}^\mathrm{MF}}
\newcommand{\Lcalo}{{\overline{\Lcal}}}

\newcommand{\Hpmf}{\mathrm{H}^\mathrm{MF}}

\newcommand{\rhoe}{{\rho_\mathrm{e}}}

\newcommand{\Ppmf}{\P^\pmf}

\newcommand{\al}{\abs{\L}}
\newcommand{\Pfrak}{\P^{\mathrm{apx}-\pmf}}
\newcommand{\EPfrak}{\E^{\mathrm{apx}-\pmf}}
\newcommand{\Zfrak}{Z^{\mathrm{apx}-\pmf}}

\newcommand{\gk}[1]{\left\{#1\right\}}
\newcommand{\ek}[1]{\left[#1\right]}
\newcommand{\rk}[1]{\left(#1\right)}
\newcommand{\oo}{\rk{1+o(1)}}
\renewcommand{\Re}{\mathrm{R}_\mathrm{e}}
\newcommand{\Oe}{\Omega_\mathrm{e}}
\newenvironment{proofsect}[1] 
{\vskip0.1cm\noindent{\scshape #1.}\hskip0.5cm} 

\begin{document}

\title[\hfill Emergence of interlacements from the finite volume Bose soup\hfill]
{Emergence of interlacements from the finite volume Bose soup}

\author[Quirin  Vogel]{Quirin  Vogel}
\address[Quirin  Vogel]{NYU Shanghai, 1555 Century Ave, Pudong, Shanghai, China, 200122\newline
Department of Mathematics, CIT, Technische Universität München, Boltzmannstr. 3, D-85748, Garching bei München, Germany.}
\email{quirinvogel@outlook.com}

\thanks{}
  

 
\begin{abstract}
    We show that conditioned on the (empirical) particle density exceeding the critical value, the finite volume Bose loop soup converges to the superposition of the Bosonic loop soup (on the whole space) and the Poisson point process of random interlacements. The intensity of the latter is given by the excess density above the critical point. We consider both the free case and the mean-field case.
\end{abstract}
 
 \maketitle
%


\section{Introduction}
The Bose gas is a physical system which has been studied both theoretically and experimentally for around 100 years. Within mathematical physics and quantum statistical mechanics, one describes the Bose gas either by employing a functional analytic framework, using \textit{CCR} (canonical commutation relation) \textit{algebras} (see \cite{bratteli2003operator}) or with the help of probabilistic models. Within those probabilistic representations, one can either work with random permutations (see \cite{ginibre1971some,betz2009spatial}) or a system of interacting random loops to compute thermodynamical properties of the Bose gas. We follow the latter approach.\\
The description of the Bose gas as a system of interacting random loops dates back to Feynman (see \cite{feynman1953atomic}); rigorous constructions can be found in \cite{ginibre1971some}. We first give an introduction into the model, while delaying the full description to Section \ref{sec_not_setup}.

\begin{figure}[h]
    \centering
    \includegraphics[scale=0.2]{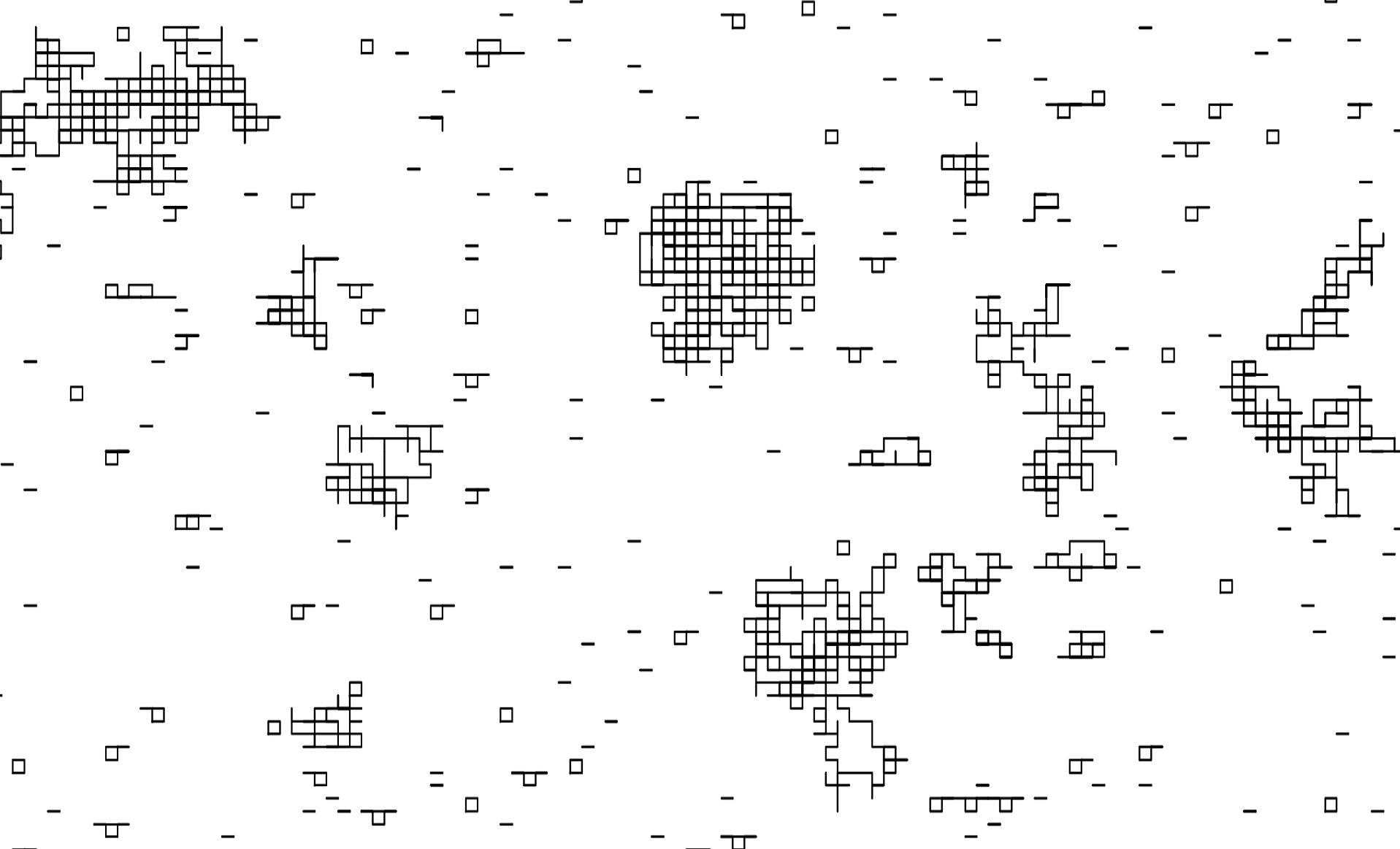}
    \caption{A sample from the loop soup on the square grid, as simulated by the author (using a Python script).}
    \label{fig:my_label1111}
\end{figure}

Fix the \textit{inverse temperature} $\beta>0$ of the system. Let $\P_{x,x}^{\beta j}$ be the unnormalized bridge measure of a continuous-time, $\Z^d$-valued simple random walk bridge from $x$ to $x$ in time $\beta j$, with $j\in\N$ and $x\in \Z^d$, see also Equation \eqref{EquationOfBridgeMeasure} for a definition. We restrict ourselves to the cases $d\ge 3$, the most relevant dimensions for the Bose gas. The \textit{Bosonic loop measure} at inverse temperature $\beta>0$ and \textit{chemical potential} $\mu\le 0$ is given by
\begin{equation}\label{FirstEquationOfM}
    M_{\L,\beta,\mu}=\sum_{x\in\L}\sum_{j\ge 1}\frac{\ex^{\beta \mu j}}{j}\P_{x,x}^{\beta j}\, ,
\end{equation}
for $\L\subset \Z^d$, on the space of loops $\gk{\omega\colon [0,\beta j]\to\Z^d\text{ with }j\in\N\text{ and }\omega(0)=\omega(\beta j)}$. The (grand canonical) model for the \textit{free} Bose gas (i.e., in the absence of interaction between particles) in the domain $\L$ is given by the Poisson point process (PPP) $\P_{\L,\beta,\mu}$ with intensity measure $M_{\L,\beta,\mu}$. This model (at times enriched with additional interactions) has been studied extensively in the literature, see for example \cite{bratteli2003operator,adams2011free,adams2018large, AV19,fröhlich2020interacting}. A sample $\eta=\sum_{\omega}\delta_\omega$ from $\P_{\L,\beta,\mu}$ is given by a collection of random loops $\omega$, also referred to as \textit{loop soup}, encoded in a point measure, see Figure \ref{fig:my_label1111} for an illustration.\\
The \textit{expected particle density} $\rho_\L(\beta,\mu)$ is the relevant macroscopic quantity for this article. For the purpose of the introduction, we define it as the expected average loop length per site, for $\L\subset\Z^d$ non-empty and finite, as
\begin{equation}\label{Equation22523}
    \rho_\L(\beta,\mu)=\E_{\L,\beta,\mu}\left[\frac{1}{\abs{\L}}\sum_{\omega\in \eta}\ell(\omega)\right]\, ,
\end{equation}
where $\ell(\omega)$ is the \textit{length} of the loop $\omega$ (if $\omega$ was sampled from $\P_{x,x}^{\beta j}$, we then set $\ell(\omega)= j$) and $\E_{\L,\beta,\mu}$ is the expectation with respect to $\P_{\L,\beta,\mu}$. We say that a loop $\omega$ with length $\ell(\omega)=j$ represents $j$ particles. From now on, abbreviate the \textit{average particle number} (per site) $\frac{1}{\abs{\L}}\sum_{\omega\in \eta}\ell(\omega)$ by $\Lcalo_\L$. In this article, we study the limit law of $\P_{\L,\beta,\mu}$ \textit{conditioned} on the event that $\{\Lcalo_\L=\rho\}$ for $\rho>0$. We do this for the free case as well as for the Bose gas with \textit{mean-field} interaction, described in the next section. We now briefly describe the main results of this paper (in the free case) and relate it to previous results. \vspace{3mm}\\
To prove a phase transition for the Bose gas, several criteria can be employed. Showing a breakdown of \textit{equivalence of ensembles} (see \cite{huang1987statistical}) or proving \textit{off-diagonal long range order} (see \cite{penrose1956bose,yang1962concept}) are two (at times equivalent) methods. For probabilistic representations of the Bose gas, a phase transition is conjectured to coincide with the existence of "infinite" loops (see \cite{feynman1953atomic,suto1,benfatto2005limit,ueltschi2006feynman,armendáriz2019gaussian}). This is usually proven via a mass-gap argument, which can be sketched as follows: after taking a certain infinite volume limit, one shows that the sum over all finite loop length probabilities is strictly less than 1. The remaining probability mass is then associated to infinite loops, see \cite{suto1,suto2,ueltschi2006feynman}.

In Theorem \ref{THM1}, we give the scaling limit of $\P_{\L,\beta,\mu}$ conditioned on the event $\{\Lcalo_\L=\rho\}$. There is a unique $\rho_\mathrm{c}>0$ such that for $\rho>\rho_\mathrm{c}$, the conditioned process becomes the independent superposition of the free Bose gas on $\Z^d$ and the process of random interlacements: for all $\mu\le 0$ and for $\rho=\rhoc+\rhoe$
\begin{equation}\label{EquationIntuitionIntro}
\lim_{\L\uparrow\Z^d}  \P_{\L,\beta,\mu}(\,\cdot\,|\Lcalo_\L=\rho_\mathrm{c}+\rho_\mathrm{e})=\P_{\Z^d,\beta,0}\otimes\P_{\rho_\mathrm{e}}^\iota\, ,  
\end{equation}
where $\P_{\rho_\mathrm{e}}^\iota$ is the Poisson point process of random interlacements with intensity $\rho_\mathrm{e}=\rho-\rho_\mathrm{c}>0$, see Section \ref{sectionInterlacements} for a definition.
The critical value $\rho_\mathrm{c}$ is equal to the limit of $\rho_\L(\beta,0)$ as $\L\uparrow\Z^d$, see Equation \eqref{Equation22523}. This is an extension of the mass-gap approach in the sense that we can explicitly describe the limiting process and see that the infinite loops coincide with the random interlacements in a path-wise sense, as conjectured in \cite{armendáriz2019gaussian}. The Markov property of the random walk also implies that $\P_{\Z^d,\beta,\mu}\otimes\P_{\rho_\mathrm{e}}^\iota$ is a Gibbs measure.

The superposition of the Bose loop soup and the random interlacements has been studied in the important work \cite{armendáriz2019gaussian}. In that work, the authors showed that superpositioning the {non-interacting} Bose loop soup \textit{on the whole space} and the random interlacements produces the same distribution on random permutations as described in \cite{macchi1975coincidence}. This suggests that the superposition of finite loops and random interlacements is the right process to study. However, in their work they did not show that this process arises from the finite volume ensemble. This was stated as an open question, resolved in this work. Furthermore, \cite{armendáriz2019gaussian} only considered the positions of the \textit{particles}, that is, the distribution of $\left\{\omega(0),\,\omega(\beta),\, \omega(2\beta),\,\ldots\right\}$. In our work, we consider the distribution of the loops/interlacements as samples in path-space. This is important, as physically relevant interactions are given via integrals of the trajectories of the loops, see \cite{ueltschi2006feynman} for example.

We also show that Equation \eqref{EquationIntuitionIntro} continues to hold under \textit{mean-field} interaction. That is, when we add interaction to the loops which depends on the square of $\Lcalo_\L$. Mean-field models tend to represent the true physical reality more accurately than free models. The mean-field model for the Bose gas is defined in Equation \eqref{EquationDefineThe MFmodel}. We refer the reader to \cite{benfatto2005limit} for a derivation of the model itself.

To our best knowledge, there is no prior work proving the emergence of random interlacements from the finite volume distributions of the Bose gas. The conditioning $\{\Lcalo_\L=\rho\}$ transfers the setting from the \textit{grand canonical} ensemble to the \textit{canonical} ensemble. For more on the physical background, we refer the reader to Section \ref{sec_dis}.

Crucial to our result is the right choice of topology, which is invariant under permutation and reparameterization of the loops. This is to be expected, from the functional-analytic description of the Bose gas.\\
There are two parts to the proof: first, one needs to show that the conditioning forces the appearance of long loops which are Poisson and asymptotically independent from the bulk of the loops. For this we tile $\Z^d$ with boxes of increasing size. Two different scales appear in the process $N^2$ and $N^{d-2}$. This can intuitively be understood through the well-known fact that a random walker intersecting a ball of radius $N$ will hit only $N^2$ points for $d\ge 3$. Hence, $N^{d-2}$ walkers are needed to generate a positive density of hitting points.\\
The second part of the proof consists in identifying the scaling limit. For this, we define an analogue of the convergence of finite-dimensional distributions for parametrization-invariant processes. While some parts for the proof are the same for the free and mean-field versions, at several places we need different arguments. Especially the fact that in the infinite volume limit, the interlacements and the normal loops decorrelate needs special attention when working with the mean-field interaction.

In our set-up of the Bose gas, the notion of boundary condition is relevant. Our proofs are for the \textit{free-boundary conditions} (not to be confused with the free Bose gas) both in the case of no interaction and mean-field interaction. This requires a slightly more complicated set-up but does use less notation in the course of the proof. For the supercritical loop representations in different boundary conditions, we refer the reader to the forthcoming \cite{koenig2023}.

The paper is organised as follows: in Section \ref{Sec_results}, we state the main results, delaying some of the more technical definitions to later. In Section \ref{sec_not_setup}, we introduce the rigorous set-up. Section \ref{sec_dis} gives the physical background and motivation in more detail; it is independent from the rest of the paper. The proof of the main results is given in Section \ref{sec_proofs} and is subdivided into the free and the interacting case.

In the Appendix, we give a table with the frequently used notation.
\section{Results}\label{Sec_results}
We begin by introducing the set-up: for $N\ge 1$, fix $\L=[-N/2,N/2)^d\cap\Z^d$, where we omit the dependency of $\L$ on $N$, to aid readability. From now on, and for the rest of the article, $\L$ will always denote this set.

Recall that $\omega\colon[0,\beta j]\to\Z^d$ denotes a loop\footnote{We later alter this definition on a set of measure zero, see Equation \eqref{MorePreciseDefininiton}.} ($\omega(0)=\omega(\beta j)$), $\ell(\omega)=j$ is the length of the loop and $\eta$ is a collection of loops encoded as $\eta=\sum_\omega\delta_\omega$. Here, $\delta_\omega$ is the delta-measure on the loop $\omega$. Define for any domain $\Delta\subset\Z^d$
\begin{equation}
     \Lcal_\Delta(\eta)=\Lcal_\Delta=\sum_{\omega\in\eta}\ell(\omega)\1\{\omega(0)\in\Delta\}
\end{equation}
Furthermore, for any domain $\Delta\subset\Z^d$, $\beta>0$, $\mu\le 0$, $\P_{\Delta,\beta,\mu}$ is the Poisson point process with intensity measure $M_{\Delta,\beta,\mu}$ defined in Equation \eqref{FirstEquationOfM}.

We now define the loop process on the \textit{whole space}: fix $\rho,\beta>0$, $\mu\le 0$ and set $\P_N^{\rho}$ to be
\begin{equation}\label{EquationDefinitionrhon}
    \P_N^{\rho}=\P_{N,\beta,\mu}^{\rho}=\bigotimes_{x\in\Z^d}\P_{xN+\L,\beta,\mu}\left(\,\,\cdot\,\,|\Lcal_{xN+\L}=\rho\al\right)\, ,
\end{equation}
where from now on we read $\Lcal_{\Delta}=x$ as $\Lcal_{\Delta}=\floor{x}$, if $x$ is not an integer, for any set $\Delta\subset\Z^d$. In words, $\P_N^{\rho}$ is the independent superposition of translates of $\P_{\L,\beta,\mu}$ conditioned to have a particle density of $\rho>0$. This gives us a distribution of loops on $\Z^d$ (rather than a finite box), allowing us to capture long range effects. We refer the reader to Remark \ref{remakr scale} for a heuristics of this set-up. Furthermore, as $\beta>0$ is fixed throughout the article and $\P_{N,\beta,\mu}^{\rho}$ does not depend on $\mu$ (as we will see in the beginning of Section \ref{sec_proofs}), we omit it from the notation and simply write $\P_{N}^{\rho}$.

Convergence in the \textit{topology of local convergence} is denoted by $\xrightarrow{\mathsf{loc}}$. It is generated by continuous and bounded functions, which depend on the values of the loops inside a compact set and are invariant under time-shifts. For a rigorous definition, we refer the reader to Definition \ref{definitionLocalconvergence}.
\begin{theorem}\label{THM1}
Fix $\beta>0$, $\mu\le 0$ and $d\ge 3$. The measure $\P_N^{\rho}$ is the law of the conditioned Bosonic loop soup on the whole space, defined in Equation \eqref{EquationDefinitionrhon}. As $N\to\infty$, it holds that if $\rho>\rhoc$
    \begin{equation}
\P_N^{\rho}\xrightarrow{\mathsf{loc}}\P_{\Z^d,\beta,0}\otimes\P_{\rho-\rho_\mathrm{c}}^\iota\, ,
    \end{equation}
where, with $\P_x$ the law of the continuous-time simple random walk on $\Z^d$,
\begin{equation}\label{EquationTHmDefinitionRhoc}
    \rho_\mathrm{c}= \rho_\mathrm{c}(\beta)=\sum_{j\ge 1}\P_x\rk{\omega(\beta j)=x}\, ,
\end{equation}
and we refer the reader to Section \ref{sectionInterlacements} for a definition of the random interlacement process $\P_{\rho-\rho_\mathrm{c}}^\iota$. The symbol $\otimes$ refers to the superposition of the two point processes, i.e., if $\eta=\sum_{j}\delta_{\omega_j}$ is sampled from $\P_{\Z^d,\beta,0}$ and $\theta=\sum_{i}\delta_{\kappa_i}$ independently from $\P_{\rho-\rho_\mathrm{c}}^\iota$, then $\eta+\theta$ is a sample from $\P_{\Z^d,\beta,0}\otimes\P_{\rho-\rho_\mathrm{c}}^\iota$.
\end{theorem}
The case $\rho<\rhoc$ is given in Theorem \ref{ThMWholeDescription}, as a variety of new objects need to be introduced first. However, no interlacements are present in the scaling limit and we consider Theorem \ref{THM1} the main result of the article.

The next result shows that the above remains true in the presence of mean-field interactions: for some scaling constant $a>0$ and $\Delta\subset\Z^d$ non-empty\footnote{we henceforth assume all sets $\L,\Delta,K\subset\Z^d$ to be non-empty.} and bounded, define $\Hpmf_\Delta$ 
\begin{equation}
    \Hpmf_\Delta(\Lcal_\Delta)=\Hpmf_\Delta(\Lcal_\Delta(\eta))=a\frac{\Lcal_\Delta^2}{2\abs{\Delta}}\, .
\end{equation}
This function is called the \textit{(particle) mean-field} Hamiltonian, see \cite{benfatto2005limit,adams2018large}. Define the corresponding probability measure $\P^\pmf_{\Delta}$ via its Radon--Nikodým derivative
\begin{equation}\label{EquationDefineThe MFmodel}
    \d \P^\pmf_{\Delta}(\eta)= \d \P^\pmf_{\Delta,\beta,\mu}(\eta)=\frac{1}{Z_\Delta^{\pmf}}\ex^{-\Hpmf_\Delta(\Lcal_\Delta(\eta))}\d \P_{\Delta,\beta,\mu}(\eta)\, ,
\end{equation}
where $Z_\Delta^{\pmf}=\E_{\Delta,\beta,\mu}\ek{\ex^{-\Hpmf_\Delta}}$.

We now define for $\beta, \rho>0$, $N\in\N$ and $\mu\le 0$
\begin{equation}\label{EquationDefPMF}
    \P_N^{\pmf,\rho}=\bigotimes_{x\in\Z^d}\P_{xN+\L,\beta,\mu}^\pmf\left(\,\,\cdot\,\,|\Lcal_{xN+\L}\ge \rho\al\right)\, ,
\end{equation}
where we omit the dependency on $\beta,\mu$ from the notation. The reason for conditioning on the $\Lcal_{xN+\L}$ being larger than $\rho\al$ is that the $ \Hpmf_\Delta$ is deterministic on the set $\gk{\Lcal_\Delta=k}$, for any $k$. Conditioning on $\Lcal_{xN+\L}$ being larger than $\rho\al$ gives similar scaling results as compared to Theorem \ref{THM1}. Formally, this corresponds to the Hamiltonian $\Hpmf_\Delta+\infty\1\gk{\Lcal_\Delta<\rho\abs{\Delta}}$ in the grand-canonical setting.
\begin{theorem}\label{THM2}
Let $\beta,\rho>0$, $\mu=0$ and $d\ge 3$. Consider the law $\P_N^{\pmf,\rho}$, defined in Equation \eqref{EquationDefPMF}. We then have that, as $N\to\infty$,
    \begin{equation}
\P_N^{\pmf,\rho}\xrightarrow{\mathsf{loc}}\P_{\Z^d,\beta,0}\otimes\P_{\rho-\rho_\mathrm{c}}^\iota \qquad\textnormal{if}\qquad\rho>\rho_\mathrm{c}\, .
    \end{equation}
\end{theorem}
\begin{remark}
The statement for the case $\rho<\rho_\mathrm{c}$ is given in \cite[Theorem 2.7]{dickson2021formation}. Interlacements do not appear in the limiting process. In \cite{dickson2021formation} the class of mean-field operators is enlarged and does incorporate all $\mu\in\R$.
\end{remark}

\section{Notation and set-up}\label{sec_not_setup}
\subsection{Domains and random walks}\label{subsectionDomain}
Fix $\L=[-N/2,N/2)^d\cap\Z^d$, the cube of sidelength $N$. Write $[n]$ to denote the set $\{1,\ldots,n\}$. 
Enrich $\Z^d$ with an additional point, denoted by the symbol $\ddagger$. The origin in $\Z^d$ is denoted by $0$.

Let $\Gamma$ be the space of doubly-infinite càdlàg (right-continuous with left limits) paths $\omega$ with values in $\Z^d\cup \{\ddagger\}$ which satisfy one of the two conditions below:
\begin{enumerate}\label{MorePreciseDefininiton}
    \item For all $t$, $\omega(t)\in\Z^d$. Furthermore, $\lim_{\abs{t}\to\infty}\abs{\omega(t)}=+\infty$.
    \item $\omega(t)\in\Z^d$ for $t\in [0,\beta j)$ and $\omega(t)=\ddagger$ for all $t\in (-\infty,0)\cup[\beta j,\infty)$ for some $\N\ni j\ge 1$.
\end{enumerate} 
Denote the subset of $\Gamma$ satisfying the first condition by $W$ and by $\Gamma_B$ for the second condition (where we write $\Gamma_{B,j}$ whenever we want to specify $j$). We think of $\Gamma_B$ as the space of finite loops (continued in the càdlàg way at $\ddagger$) and $W$ the space of doubly infinite paths (a model for our interlacements). For a loop $\omega\in\Gamma_{B,j}$, we denote its \textit{length} by $\ell=\ell(\omega)= j$. The physical interpretation is that such a loop \textit{represents} $j$ particles, see Section \ref{sec_dis}. We work with the Skorokhod topology on $\Gamma$ and refer the reader to \cite[Section 16]{billingsley1968convergence} for an explicit construction of the metric inducing this topology. Our sigma-algebra is the Borel sigma-algebra.

Given $t\in \R$ we define the shift by $t$, denoted by $\theta_t$, as follows: for doubly-infinite paths, $\theta_t$ acts as a normal shift and for finite loops, we shift the trajectory but keep the loop parameterised on $[0,\beta j)$, i.e.,
\begin{enumerate}
    \item $\omega\circ\theta_t(s)=\omega(t+s)$, if $\omega\in W$.
    \item $\omega\circ\theta_t(s)=\omega(t+s\mod \beta j)$, if $\omega$ $\in$ $\Gamma_{B,j}$.
\end{enumerate}
Denote $\Gamma^*$ to be $\Gamma$ modulo the shift operator, i.e., the space of equivalence classes of $\Gamma$, where $\omega_1$ is equivalent to $\omega_2$ if there exists $t$ with $\omega_1=\omega_2\circ\theta_t$. In the same way we set $W^*$ and $\Gamma_{B}^*$ the respective subspaces of equivalence classes. Let $\Pi$ denote the projection from $\Gamma$ to $\Gamma^*$ and let $\Pi^{-1}[A]=\gk{\omega\in\Gamma\colon\Pi(\omega)\in A}$ for any set $A\subset \Gamma^*$.

For $x\in\Z^d$, recall that $\P_x$ is the law of the continuous-time, nearest-neighbor simple random walk started at $x$. We denote its transition kernel by $p_t(x,y)=p_t(x-y)$, for $t>0$ and $x,y\in\Z^d$. For $x\in\Z^d$ and $j\in \N$, define the bridge measure $\P_{x,x}^{\beta j}$ via 
\begin{equation}\label{EquationOfBridgeMeasure}
    \d\P_{x,x}^{\beta j}=\1\left\{\lim_{s\uparrow {\beta j}}\omega(s)=x\right\}\d \P_x\otimes\1\left\{\omega(s)=\ddagger,\,\forall s\in\R\setminus [0,\beta j)\right\}\, .
\end{equation}
In words, for $s\in [0,\beta j)$, the law of $\omega(s)$ is governed by the unnormalized bridge measure induced by the continuous-time simple random walk from $x$ to $x$ in time $\beta j$, for all other $s$, we set $\omega(s)=\ddagger$. We also define $\B_{x,x}^t=p_t(x,x)^{-1}\P_{x,x}^t$, the normalised bridge measure. We work with the interval $[0,\beta j)$ instead of $[0,\beta j ]$ so that the resulting path remains càdlàg. As the probability of a jump at time $\beta j$ is zero, this does not affect the transition probabilities. The set-up is such that we can formally have interlacements and loops living on the same space, and the distributions converging to one another.

For a set $K\subset \Z^d$, let $H_K$ denote the first time (after its first jump) the random walk hits $K$
\begin{equation}
    H_K=\inf\{t>t_o\colon \omega(t)\in K\}\, ,\quad\textnormal{where}\quad t_o=\inf\{t>0\colon \omega(0)\neq \omega(t)\}\, .
\end{equation}
We will write $H_x$ instead of $H_{\gk{x}}$ throughout the article, for $x\in\Z^d$.

We furthermore write $\Rcal_t$ for the number of distinct vertices the random walker has visited up to time $t\ge 0$
\begin{equation}
    \Rcal_t=\Rcal_t(\omega)=\#\gk{x\in\Z^d\colon \exists 0<s\le t\text{ with }\omega(s)=x}\, .
\end{equation}

For $\Delta\subset\Z^d$ (recall that all sets $\L,\Delta, K$ are assumed to be non-empty), $\beta>0$, and $\mu\le 0$, recall the definition of the Bosonic loop measure from Equation \eqref{FirstEquationOfM}.

Most of the time, we will work with $\mu=0$. As $\beta>0$ remains fixed throughout this article, we do not incorporate it into the notation. Furthermore, if $\mu$ is omitted as index, it means that it is zero. This means that we write $M_\Delta$ for $M_{\Delta,\beta,0}$. Note that for $M_\Delta$, the random walk bridges start within $\Delta$ but are allowed to exit $\Delta$. This will be important later on. 

Let $\P_{\Delta,\beta,\mu}$ be the law of the Poisson point process (PPP) with intensity measure $M_{\Delta,\beta,
mu}$ for $\Delta\subset\Z^d$. Recall that $\P_\Delta$ is short for $\P_{\Delta,\beta, 0}$. A sample from $\P_{\Delta,\beta,\mu}$ is denoted by $\eta$ and can be written as
\begin{equation}
    \eta=\sum_{k}\delta_{\omega_k}\, ,
\end{equation}
where $(\omega_k)_k$ are (almost surely) distinct loops (due to the continuous-time) and we sum over finitely many $k$ if and only if $\Delta$ contains finitely many points. We refer to the collection $(\omega_k)_k$ as the \textit{loop soup}. To make notation more readable, we write $\omega\in\eta$ instead of $\omega\in\supp(\eta)$. 

As, for $d\ge 3$, the random interlacements induce two "natural scales", we need to introduce some additional notation. We introduce a domain $\L_N$ which is finite but big enough to capture the long-range effects.
We briefly comment on the scaling used in this article.
\begin{remark}\label{remakr scale}
For the free boundary condition to exceed a particle density of $\rho_\mathrm{c}$ in $\L$, the system will produce one big loop. This loop has duration of the same order as the volume of $\L$ (i.e., $N^d$) and thus its diameter is of order $N^{d/2}$. This means it stretches across a large box containing $N^{d/2-1}$ cubes of sidelength $N$. So to properly capture all the long loops intersecting the box $\L$, one needs to consider at least $N^{d^2/2-d}$ copies of the box $\L$. 


This scaling is to be expected: a typical interlacement (a single doubly-infinite random walk trajectory) intersecting the box $\L$ (or $[-N^{d/2},N^{d/2}]^d\cap\Z^d$) will visit of order $N^2$ (resp. $N^d$) points in it. Thus, we typically have to have $N^{d-2}$ (resp. $N^{d^2/2-d}$) interlacements hitting the box $\L$ to produce a non-vanishing local density.
\end{remark}
Recall that $\L=[-N/2,N/2)^d\cap\Z^d$. We choose $r_N=\log^2(N+1)$, positive, increasing and diverging to infinity. The reason for this sequence is that the long loops intersecting the origin have a diameter of around $N^{d/2}$ on average, see Remark \ref{remakr scale}. However, we consider loops started at a distance of at most $r_NN^{d/2}$ away from the origin, in order to incorporate all long loops with very high probability. Set $C_N=r_N[-N^{d/2-1}/2,N^{d/2-1}/2)^d\cap\Z^d$ (the set of centres of our cubes) and $\L_N=\bigcup_{x\in C_N}\left(xN+\L\right)$. Thus, $\L_N=r_N\left[-N^{d/2}/2,N^{d/2}/2\right)^d\cap \Z^d$. We have that $\P_{\L_N}=\P_{\L_N,\beta,0}$ is our standard loop process on the domain $\L_N$. It will be helpful to think of $\P_{\L_N}$ as the independent superposition of the laws of $\P_{xN+\L}$ with $x\in C_N$. It will be shown that the measure $\P_{\L_N}$ samples all the loops contributing to Theorem \ref{THM1} and Theorem \ref{THM2}.

Let the {occupation field} at $x\in\Z^d$ be defined as
\begin{equation}
    \Lcal_x=\Lcal_x(\eta)=\sum_{\mycom{\omega\in \eta}{\omega(0)=x}}\ell(\omega)=\eta[\ell(\omega)\1\{\omega(0)=x\}]\, .
\end{equation}
We will denote for $\Delta\subset\Z^d$
\begin{equation}\label{DefinitionLocalTime}
   \Lcal_\Delta=\sum_{x\in\Delta}\Lcal_x\qquad\text{and}\qquad \overline{\Lcal}_\Delta=\frac{1}{\abs{\Delta}}\sum_{x\in\Delta}\Lcal_x\, ,
\end{equation}
where $\Delta$ needs to be finite for the second equality to hold. Owing to the identification of $\ell(\omega)$ with the number of particles $\omega$ represents, we can think of $\Lcal_\Delta$ as the total number of particles (started) in $\Delta$ and of $\Lcalo_\Delta$ as the particle density.

Note that by the Campbell formula (see \cite[Proposition 2.7]{last2017lectures}) and the translation invariance, we have that for all $x\in\Delta$
\begin{equation}\label{MeanLocalTimeEquation}
    \E_{\Delta,\beta,0}[\Lcal_x]=\frac{1}{\abs{\Delta}}\E_{\Delta,\beta,0}\ek{\Lcal_\Delta}=\frac{1}{\abs{\Delta}}M_{\Delta,\beta,0}\ek{\ell(\omega)}=\sum_{j\ge 1}\P_x\rk{\omega(\beta j)=0}=\sum_{j\ge 1}p_{\beta j}(0)=\rhoc\, ,
\end{equation}
also compare with Equation \eqref{EquationTHmDefinitionRhoc}.

We now define the topology of local convergence.
\begin{definition}\label{definitionLocalconvergence}
We say that a function $F$ from the space of point measures on $\Gamma$ to $\R$ is \textit{path-local and symmetric} if there exists a bounded set $K\subset \Z^d$ such that:
\begin{enumerate}
    \item For all configurations $\eta$, there exist $(F_n)_n$ functions, $F_n\colon\Gamma^{\otimes n}\to \R$, such that for all $n\ge 0$
    \begin{equation}\label{reallytemp}
        F(\eta)=F_n(\omega_1,\ldots,\omega_n)\, ,
    \end{equation}
    if $\{\omega\in\eta\colon \omega\cap K\neq\emptyset\}=\{\omega_1,\ldots,\omega_n\}$.
    \item For each $n$, $F_n$ is invariant under permutation of its arguments (and hence \eqref{reallytemp} is well-defined).
    \item For all $n\in\N$, if $(\omega_i)_{i=1}^n$ and $(\tilde\omega_i)_{i=1}^n$ agree inside $K$, we have that
    \begin{equation}
 F_n\left(\left(\tilde\omega_i\right)_{i=1}^n\right)=F_n\left(\left(\omega_{i}\right)_{i=1}^n\right)\, .
\end{equation}
\end{enumerate}
We say that $F=(F_n)_n$ is \textit{shift-invariant}, if for any $n$ and $F_n\colon \Gamma^{\otimes n}\to\R$ and for any collection of $t_i\in\R$, we have that
\begin{equation}
     F_n\left(\left(\omega_i\right)_{i=1}^n\right)=F_n\left(\left(\omega_i\circ \theta_{t_i}\right)_{i=1}^n\right)\, .
\end{equation}
The \textit{topology of local convergence} on point measures on $\Gamma$ is generated by the convergence of continuous and bounded functions which are also path-local, symmetric and shift-invariant. We call such functions \textnormal{test functions}. We denote convergence in this topology by the symbol $\xrightarrow{\mathsf{loc}}$.

\end{definition}
Note that the topology is quite natural for the Bose gas: the shift-invariance and the symmetry come from the quantum mechanical description (see \cite{bratteli2003operator}). 

We give some examples of path-local functions.
\begin{enumerate}
\item Let $F(\eta)=\1\{\eta\in A\}$ with $A$ a local event, which does not depend on parametrization or order of the loops. Important cases of such $A$ include: $A=\{\textnormal{ no loops intersect }\Delta\}$, $A=\{
\text{ the origin is connected to }\Delta^c\text{ through loops }\}$, $A=\{
\text{ total local time at the origin  }\in [a,b]\}$ ($a,b\in\R$), with $\Delta\subset\Z^d$ bounded.
\item The local mean-field energy: for $\Delta\subset\Z^d$ bounded and $a>0$,
\begin{equation}
    F(\eta)=\exp\left(-\frac{a}{2}\Lcalo^2_\Delta\right)\, ,
\end{equation}
is path-local. The same holds for the energy of the HYL-model, see \cite{dickson2021formation} for a definition.
\item The modified\footnote{We say modified because the structure of $\widetilde{V}_K$ is slightly different from the interaction defined in \cite{bratteli2003operator}.} path-interaction energy inside $K$: fix $K\subset \Z^d$ and define
\begin{multline}\label{interactionK}
 \widetilde{V}_K(\omega^{\ssup{1}},\ldots,\omega^{\ssup{n}})= \\
\frac{1}{2}\sum_{1\le i<j\le n}\int_{D(\omega^{(i)})\times D(\omega^{(j)})}\,v(|\omega^{\ssup{i}}(t_i)-\omega^{({j})}(t_j)|)\1\{\omega^{\ssup{i}}(t_i),\omega^{\ssup{j}}(t_j)\in K\}\,\d t_i\, \d t_j\, ,
\end{multline}
where $D(\omega^{(j)})$ is the duration of the $j$-th path, either $[0,\beta\ell(\omega))$ if the path is a loop, or $(-\infty,\infty)$ if the path is an interlacement. For a strongly repelling interaction, i.e., $v>c>0$, we have that
\begin{equation}
    F(\eta)=F\left(\sum_{\omega}\delta_\omega\right)=\ex^{-\widetilde{V}_K\left(\eta\right)}\, ,
\end{equation}
is a path-local function.
\end{enumerate}
If $f(x)$ and $g(x)$ are two functions, we write $f(x)\sim g(x)$ (or $f\sim g$) whenever $f(x)=g(x)(1+o(1))$ in the Landau notation and the limit is apparent from the context. 

The convention used for constants is as follows: by $C, c$ we denote positive
constants depending only on the dimension $d$, whose value may change from place to place. Dependence of constants on additional parameters will be highlighted explicitly.
\subsection{Random Interlacements}\label{sectionInterlacements}
The theory of random interlacements was proposed by Sznitman in \cite{sznitman2010vacant} and has been a subject of intense study, see \cite{drewitz2014introduction} for an introduction and overview. We now give a very brief definition of the underlying process:\\
Recall that the subset $W$ of $\Gamma$ consists of those paths whose absolute value diverges to infinity as the absolute value of the time diverges to infinity
\begin{equation}
    W=\left\{\omega \in \Gamma\colon \lim_{\abs{t}\to\infty}\abs{\omega(t)}=+\infty\right\}\, .
\end{equation}
Recall that $W^*$ is equal to $W$ modulo the shift operator. Let $W_K$ (resp. $W_K^*$) be the set of those paths in $W$ (resp. $W^*$) which intersect $K$, for $K\subset\Z^d$ a bounded set. Set $\kappa_d=\P_0\rk{H_0=\infty}$, the escape probability of the random walk. Let $\mathtt{Q}^K$ be that measure on $W_K$ which has finite dimensional distributions given by
\begin{multline}\label{EquationDefInterl}
    \mathtt{Q}^K\left(\omega({-s_i})=y_i,\, \omega(0)=z,\, \omega({t_j})=x_j,\,\forall \, i\in [m],\, j\in [n]\right)\\
    =\kappa_d\P_z\left(\omega(-{s_i})=y_i,\,\forall \, i\in [m]|H_K=+\infty\right)e_K(z)\P_z\left( \omega({t_j})=x_j,\, \forall j\in [n]\right)
\end{multline}
with $-\infty<s_m<\ldots<s_1<0<t_1<\ldots<t_n<\infty$ and a set of points $y_m,\ldots,y_1,z,x_1,\ldots x_n\in\Z^d$, with $z\in K$. Here, $e_K(z)=\P_z(H_K=+\infty)$ is the (unnormalized) equilibrium measure at $z$.

Let $\nu$ be the measure on $W^*$ which satisfies for all $K\subset\Z^d$ the relation
\begin{equation}
    \nu(A)=\mathtt{Q}^K\left[\Pi^{-1}\ek{A}\right]\, ,
\end{equation}
for $A\subset W^*_K$ measurable. For the construction of such a measure, we refer the reader to \cite{drewitz2014introduction}.

For $u\ge 0$, we denote the PPP of the (continuous-time) random interlacements at level $u$ by $\P_u^\iota$, i.e., $\P_u^\iota$ is the PPP with intensity measure $u\nu$ on $W^*$. A sample from $\P_u^\iota$ can be denoted by
\begin{equation}
    \sum_{k}\delta_{\omega^*_k}\, ,
\end{equation}
with $\omega^*\in W^*$. Note (Equation \eqref{EquationDefInterl}) that our definition of random interlacements at level $u$ is equivalent to the standard definition of random interlacements at level $u\kappa_d$. We introduce this multiplicative factor so that our version of the random interlacements at level $u$ have an average density of $u$ instead of $u/\kappa_d$ as with standard definition, compare \cite[Equation 2.21]{sznitman2012}. 

In \cite[Theorem 4.14]{Sznitman} one has the following approximation result: take the Markovian loop measure $M^{\texttt{Mark}}$ (introduced in \cite{le2010markov})
\begin{equation}
    M^{\texttt{Mark}}=\sum_{x}\int_0^\infty \frac{1}{t}\P_{x,x}^t\, \d t\, ,
\end{equation}
and denote $\P_a^{\texttt{Mark}}$ the PPP point process with intensity measure $aM^{\texttt{Mark}}$, $a>0$. When, choosing $a=u R^{d-2}c_d^{-1}$, we have that
\begin{equation}
    \lim_{R\to\infty}\P_a^{\texttt{Mark}}
    \left(\exists \omega\in\eta\colon \omega\cap K\neq\emptyset\textnormal{ and }\omega\cap\bB_R^c\neq\emptyset\right)=\ex^{-u\mathsf{Cap}(K)}=\P^\iota_u(\exists\omega\in\eta\colon\omega\cap K\neq \emptyset)\, ,
\end{equation}
where $K\subset\Z^d$ and $c_d=d\Gamma(d/2-1)\pi^{-d/2}/2$ and $\mathsf{Cap}$ is the capacity of the set $K$. Theorem \ref{THM1} can be seen as a variation of the above to a path-wise level. Indeed, the result from \cite{Sznitman} is on the level of \textit{sites} visited by any loop, whereas our topology is finer. However, this result considers diverging intensity, whereas our theorems are concerned with supercritical densities.
\subsection{Some results on Poisson point processes}
In this section, we briefly recall some results on Poisson point processes. Let $\P$ be a Poisson point process with intensity measure $\nu$ on some locally compact Polish space $\rk{\Omega, \mathcal{E}}$. We assume that for any $\omega\in\Omega$, we have $\nu\rk{\gk{\omega}}=0$. We furthermore represent a sample from $\P$ as $\sum_\omega\delta_\omega$.
\begin{proposition}\label{PropositionPoisson}
Suppose that $\nu(\Omega)<\infty$ and that $\P$ is a Poisson point processes with intensity measure $\nu$ and expectation $\E$. Write $\nu[f]$ for $\int f\d \nu$, $f$ measurable.
\begin{enumerate}
    \item For any function $f\ge 0$
    \begin{equation}
        \E\ek{f(\eta)}=\ex^{-\nu\rk{\Omega}}f(0)+\ex^{-\nu\rk{\Omega}}\sum_{k\ge 1}\frac{1}{k!}\int f\rk{\delta_{\omega_1}+\ldots+\delta_{\omega_k}}\d\nu^{\otimes k}\rk{\omega_1,\ldots,\omega_k}\, .
    \end{equation}
    \item For $h\colon\Omega\to \R$ bounded, the \textnormal{Campbell formula} gives
    \begin{equation}
        \E\ek{\eta[h]}=\nu\ek{h}\qquad\textnormal{and}\qquad \E\ek{\ex^{\eta[h]}}=\ex^{\nu\ek{\ex^{h}-1}}\, .
    \end{equation}
    \item Suppose $A$ is an increasing event (i.e., if $\supp(\eta)\subset\supp(\eta')$ then $\eta\in A$ implies $\eta' \in A$) and $B$ a decreasing event (i.e., $B^c$ an increasing event), then
    \begin{equation}
        \P\rk{A\cap B}\le \P\rk{A}\P\rk{B}\, .
    \end{equation}
    \item For $A\in\Ecal$ and $B$ an event, we have
    \begin{equation}\label{Equation244}
        \P\rk{\exists\omega_o\colon \omega_o\in A,\, \eta-\delta_{\omega_o}\in B}\le \nu(A)\P\rk{\eta\in B}\, .
    \end{equation}
    Furthermore, if $\exists!$ denotes the quantifier that there exists a unique element, then by the Poisson property
    \begin{equation}\label{Equation244!}
        \P\rk{\exists !\omega_o\colon \omega_o\in A,\, \eta-\delta_{\omega_o}\in B}= \nu(A)\P\rk{\eta\in B,\forall \omega\in\eta:\omega\notin A}\, .
    \end{equation}
\end{enumerate}
\end{proposition}
The first statement follows directly from the Poisson property, see \cite[Exercise 3.7]{last2017lectures}. For the second statement, see \cite[Exercise 3.4]{last2017lectures}. For the third statement, see \cite[Theorem 20.4]{last2017lectures}. Equation \eqref{Equation244} follows from the Mecke equation, see \cite[Theorem 4.1]{last2017lectures} and a union bound. Equation \eqref{Equation244!} is the Mecke equation (applied to the function $f(x,\eta)=\1\gk{x\in A}\1\gk{\rk{\eta-\delta_x}\ek{A}=0,\, \rk{\eta-\delta_x\in B}}$ in the notation used in the reference).

\section{Discussion}\label{sec_dis}
\subsection{Bose--Einstein Condensation}\label{sub_sec_BEC}
We do not give a full account of Bose--Einstein condensation (BEC) but instead, provide some motivation for our main theorem. For a rigorous in-depth account, we refer the reader to \cite{bratteli2003operator}.

Bosons constitute one type of fundamental particles, with fermions being their counterpart. Examples of bosons are photons, gluons or the Higgs boson. However, composite particles, such as atoms, can also behave like bosons. Bosons do not obey the Pauli exclusion principle and tend to aggregate. Bose--Einstein condensate is a (physical) state of matter in which a macroscopic fraction of particles occupy the same quantum state. It occurs at sufficiently high densities (or low temperatures). Bose--Einstein condensation is a quantum phase transition, which had been predicted by Bose and Einstein (see \cite{bose1924plancks, einstein2005quantentheorie}) decades before it was observed in a laboratory (see \cite{anderson1995observation,davis1995bose}), due to the extreme conditions needed for it to appear. For the free Bose gas at density $\rho>0$, the calculations of Bose and Einstein predict this fraction to be $\max\{\rho-\rho_\mathrm{c},0\}$.\\
BEC remains an active topic of investigation to experimental physicists, theoretical physicists and mathematicians alike.

Several criteria to prove BEC from a mathematical model have been proposed. Amongst those, \textit{off-diagonal long-range order} (ODLRO) (see \cite{penrose1956bose,huang1957quantum}) and \textit{breakdown of equivalence of ensembles} (see \cite{huang1987statistical}) are prominent. ODLRO means that the 1-particle-reduced density matrix has an eigenvalue of order volume while breakdown of equivalence of ensembles is the non-equivalence (on the level of thermodynamical functions or on the level of ensembles) of the grand canonical and the canonical ensemble after taking the thermodynamic limit. Ever since Feynman brought forward his stochastic representation of the Bose gas, it has been conjectured that (at least in many physically relevant cases), BEC coincides with an occurrence of large loops (see \cite{feynman1953atomic,van1990pressure,benfatto2005limit,ueltschi2006feynman,armendáriz2019gaussian} for example). We briefly outlined past arguments in the introduction. We should also mention that the occurrence of the Poisson--Dirichlet distribution is taken as an indicator of BEC, see \cite{betz2009spatial} for example.

In this work, we do \textit{not} prove the equivalence of different criteria for BEC. This is subject to future work. We merely calculate the limit of the finite volume representations of the Bose gas. The overarching assumption in statistical mechanics is that after taking the thermodynamic limit, the resulting process represents the physical reality on a qualitative level. The fact that our limiting process puts positive weight on random interlacements for densities above the critical point gives heuristic evidence to support the conjecture that "infinite loops" indeed represent the condensate fraction of the gas. Furthermore, this positive weight agrees with the density calculated by Bose and Einstein. We would like to point the reader to the work \cite{dickson2021formation}, in which it is shown that for a non mean-field interaction, random interlacements appear in the thermodynamic limit, thus supporting the assumption that this behavior can be expected for a wide range of interactions, see Section \ref{sectionInt} for recent results in that direction.
\subsection{Lattice vs continuum}
We chose to use the lattice for our computations. However, our proof does not depend on that: one can define with some effort, the Brownian interlacements, see \cite{sznitman2013scaling}. This is demonstrated in \cite{dickson2021formation}, where the results of this paper are generalised to different interactions, for the Brownian loop soup.
\subsection{The chemical potential}\label{subsectionChem}
For $\mu\le 0$, we introduce the following (thermodynamic) functions: the density
\begin{equation}\label{Definitionbrho}
    \boldsymbol{\rho}(\mu)=\sum_{j\ge 1}{\ex^{\beta \mu j}}p_{\beta j}(0)\, ,
\end{equation}
and the pressure
\begin{equation}\label{EquationIntroductionPressure}
   \Phi(\mu)=\sum_{j\ge 1}\frac{\ex^{\beta \mu j}}{j}p_{\beta j}(0)\, .
\end{equation}
Note that $\rho_\mathrm{c}=\boldsymbol{\rho}(0)$. Furthermore, for $0<x\le \rho_\mathrm{c}$, we set $\boldsymbol{\mu}(x)$ to be the unique $\mu\le 0$, such that $\rho(\mu)=x$. Thus, $\boldsymbol{\rho}(\boldsymbol{\mu}(x))=x$.
We now give the results for the case $\mu<0$.
\begin{theorem}\label{ThMWholeDescription}
Recall that $\P_N^{\rho}$ was defined in Equation \eqref{EquationDefinitionrhon} and is the law of $\P_{N,\beta,\mu}$ conditioned on $\{\Lcalo_\Delta=\rho\}$ in each box $\Delta=xN+\L$. We then have for $\mu\le 0$, $\beta>0$, $d\ge 3$ and $\rho<\rhoc$
\begin{equation}
    \P_N^\rho\xrightarrow{\mathsf{loc}}\P_{\Z^d,\beta,\boldsymbol{\mu}(\rho)}\qquad\textnormal{as }N\to\infty \, .
\end{equation}
\end{theorem}
The proof of the theorem above is given in \cite[Theorem 2.7]{dickson2021formation}. In the interest of being self-contained, we give it in the appendix.

The above together with Theorem \ref{THM1} may be summarized colloquially as follows:
\begin{center}
    \textit{For the loop ensemble to achieve a large particle density, it will first reduce the chemical potential (up to zero). If that cannot provide sufficient additional density, it will start to produce random interlacements.}
\end{center}
Theorem \ref{ThMWholeDescription} also indicates a break-down of equivalence of ensembles: for densities less than $\rho_\mathrm{c}$, the limiting measure of the canonical ensemble is given by the grand canonical ensemble with adjusted chemical potential $\mu=\boldsymbol{\mu}(\rho)$. For $\rho>\rho_\mathrm{c}$, this breaks down as the limiting process involves random interlacements.
\begin{remark}
Note that our results allow for a possible formulation of the infinite loops criterion:
\begin{center}
   Fix $\beta>0$. Conditioned on the empirical density $\Lcalo_\L$ being approximately equal to $\rho$, do we observe interlacements in the limit?\\ 
\end{center}
\end{remark}
\subsection{Interactions}\label{sectionInt} In this work, we show that the emergence of interlacements holds for both the free case and the mean-field case. While some parts of the proof for the free case carrry over, the behavior (and thus the proof) in the mean-field case is different. We expect universality however: for a large class of interactions, the same phenomena observed in Theorem \ref{THM1} and Theorem \ref{THM2} should persist, see \cite{dickson2021formation} for a non mean-field model.\\
However, given the presence of interaction terms, estimating the contribution of the long loops becomes increasingly complex the more the interaction depends on the geometry of the random path. In a future publication, we will show that Theorem \ref{THM2} remains valid for a certain (restrictive) class of Hamiltonians. The recent works \cite{Deuchert_2021,collin2022free,quitmann2022macroscopic} seem to have made some progress in that direction.
\section{Proofs}\label{sec_proofs}
We will only prove the result for the free boundary conditions, as mentioned in the introduction.


From now on, we assume $\rho=\rho_\mathrm{c}+\rho_\mathrm{e}$ with $\rho_\mathrm{e}>0$.
    
The proof is structured as follows:
\begin{enumerate}
    \item We first recall some results on heavy-tailed random variables from the literature. We then transfer them to our setting.
    \item Next, we prove a number of technical results, essentially ensuring that the long loops hitting the origin behave typically with high probability.
    \item We then show that the loops are asymptotically governed by a Poisson point process. We then show that the intensity measure of that process converges to that of random interlacements.
\end{enumerate}
Before proceeding further, we show that we can assume $\mu=0$ without loss of generality: note that
\begin{equation}
    \frac{\d\P_{\L,\beta,\mu}}{\d\P_{\L,\beta,0}}(\eta)=\ex^{-\al\ek{\Phi(\mu)-\Phi(0)}}\ex^{\beta\mu\Lcal_\L(\eta)}\, ,
\end{equation}
using the first statement of Proposition \ref{PropositionPoisson} and Equation \eqref{EquationIntroductionPressure}. Hence, for any test function $F$
\begin{equation}\label{EquationChangeOfMeasure}
    \E_{\Lambda,\beta,\mu}\ek{F|\Lcal_\L=\rho\al}=\frac{\E_{\L,\beta,0}\ek{F,\Lcal_\L=\rho\al}}{\P_{\L,\beta,0}\rk{\Lcal_\L=\rho\al}}=\E_{\Lambda,\beta,0}\ek{F|\Lcal_\L=\rho\al}\, ,
\end{equation}
where we used that on the event $\gk{\Lcal_\L=\rho\al}$, the Radon--Nikodým derivative becomes constant. This was to be expected as the conditioning takes us from the grand-canonical ensemble (indexed by $\beta$ and $\mu$) to the canonical ensemble (indexed by $\beta$ and $\rho$). We henceforth assume that $\mu=0$ unless explicitly stated otherwise.
\subsection{Some preliminary technical results}
In this section, we give some background on the behaviour of sums of heavy-tailed random variables. Our goal is to derive precise estimates of the probability that $\Lcalo_\L$ is much larger than its mean and gain insight into how this event is achieved.
\begin{lemma}\label{LemmaRefCLT}
Suppose $\rk{X_i}_i$ is a collection of independent, identically distributed, lattice-valued\footnote{i.e., $(X_i)_i$ takes values in $\gk{ax+h\colon x\in\Z}$ for some $a,h>0$ real numbers.} random variables such that for some $d\in\N$ with $d\ge 3$ and some $p>0$
\begin{equation}\label{EquationAssumptionXi}
    \P\rk{X_1=x}=\frac{p}{x^{1+d/2}}\rk{1+o(1)}\, ,
\end{equation}
as $x\to\infty$ and also $\P\rk{X_1<-M}=0$ for $M>0$ large enough. Then $\rk{X_i}_i$ satisfies a central limit theorem at scale $\rk{a_n}_n$ with 
\begin{equation}
    a_n=\begin{cases}
    n^{2/3}&\textnormal{ if }d=3\, ,\\
    \sqrt{n\log n}&\textnormal{ if }d=4\, ,\\
    n^{1/2}&\textnormal{ if }d\ge 5\, ,\\
    \end{cases}
\end{equation}
where for $d=3$ the limit is the centered $\alpha$-stable distribution with $\alpha=3/2$ and for $d\ge 4$ the limit is Gaussian with some finite, dimension-dependent variance.
\end{lemma}
\begin{proof}
The above result is standard, see for example \cite[Theorem 8.3.1]{bingham1989regular} or \cite[Theorem 16.29]{klenke2013probability} for a reference: to verify the conditions in these reference, we need to calculate the truncated second moment
\begin{equation}
    U(y)=\E\ek{X_1^2\1\gk{\abs{X_1}<y}}\, .
\end{equation}
Equation \ref{EquationAssumptionXi} shows that $U(y)$ grows (up to some multiplicative constant) like $y^{1/2}$ for $d=3$, $\log(y)$ for $d=4$ and is uniformly bounded for $d\ge 5$. Hence, for the case $d=3$, $U(y)y^{3/2-2}$ is slowly varying (in this case constant), which proves the result. For $d\ge 4$, note that $U(y)$ is slowly varying, and hence the limit is Gaussian. This concludes the proof. 
\end{proof}
 Next, we collect some results on the probability that the sum of such $\rk{X_i}_i$ is large. The next lemma states that if the sum of $\rk{X_i}_i$ is large, there exists a single $X_i$ realising this large value. 
\begin{figure}
    \centering
    \begin{tikzpicture}[very thick]
    \draw[black] (0,0) -- (8,0) node[below]{$\Z$};
    \draw [decorate,decoration = {calligraphic brace,amplitude=5pt}] (1,0.25) --  (2,0.25)node[pos=0.5,above=5pt,black]{\small{$S_{n-1}$}};;
    \draw [decorate,decoration = {calligraphic brace, mirror,amplitude=2pt}] (1,-0.5) --  (1.5,-0.5)node[pos=0.5,below=2pt,black]{\small{$a_n^{1+\delta}$}};
    \draw[thin] (1.5,-0.075)node[below]{\small{$0$}}  -- (1.5,0.075);
    \draw[thin] (2,-0.075) -- (2,0.075);
    \draw[thin] (1,-0.075) -- (1,0.075);
    \draw [decorate,decoration = {calligraphic brace,amplitude=5pt}] (6,0.25) --  (7,0.25)node[pos=0.5,above=5pt,black]{\small{$M_n$}};;
    \draw[thin] (6.5,-0.075)node[below]{\small{$x$}}  -- (6.5,0.075);
    \draw [decorate,decoration = {calligraphic brace, mirror,amplitude=2pt}] (6,-0.5) --  (6.5,-0.5)node[pos=0.5,below=2pt,black]{\small{$a_n^{1+\delta}$}};
    \draw[thin] (7,-0.075) -- (7,0.075);
    \draw[thin] (6,-0.075) -- (6,0.075);
    \end{tikzpicture}
    \caption{A sketch of the behaviour of heavy-tailed random variable, conditioned on the sum being large ($S_n=x$):  the maximum $M_n=\max_{i\in[n]}X_i$ will be close to $x$, with an error bounded by $a_n^{1+\delta}$ with high probability. The remaining sum will fluctuate around its mean.}
    \label{fig:my_label1}
\end{figure}
The fluctuations of that $X_i$ are controlled by the underlying CLT. See Figure \ref{fig:my_label1} for an illustration. Recall that $[n]=\gk{1,\ldots,n}$.
\begin{proposition}\label{PropositionLDP}
Let $\rk{X_i}_i$ be as in the previous lemma and assume that $\E\ek{X_1}=0$. Fix $d\ge 3$. Let $S_n=\sum_{i=1}^nX_i$. Fix $1>c>0$. There exists $\delta>0$ such that for all $x$ with  $c^{-1}n\ge x\ge c n$
\begin{equation}\label{EquationFirst1}
    \P\rk{S_n=x}\sim n\P\rk{X_1=x}\sim p n x^{-1-d/2}\qquad\text{as }n\to\infty\, ,
\end{equation}
as well as
\begin{equation}\label{EquationFirst2}
    \P\rk{S_n=x, \nexists i\in [n]\colon x-a_n^{1+\delta}\le X_i\le x+a_n^{1+\delta}}=o(1) \P\rk{S_n=x}=o\rk{n^{-d/2}}\, .
\end{equation}
We stress that the $1+o(1)$ implicit in the two equations above, does not depend on $x$. Furthermore, there exists $C>0$ such that for all $x\ge a_n^{1+\delta}$
\begin{equation}\label{EquationUnifromBound}
    \P\rk{S_n=x}\le C n\P\rk{X_1=x}\, .
\end{equation}
\end{proposition}
\begin{proofsect}{\textbf{Proof of Proposition \ref{PropositionLDP}}}
All proofs in the references are based on the Fuk--Nagaev inequality, which controls the probability of the sum being large but the maximum being small, see \cite{fuk1971probability,nagaev1982asymptotic}. For Equation \eqref{EquationFirst1} and Equation \eqref{EquationUnifromBound}, see \cite[Theorem 2.4]{berger2019notes} for the case $d=3$, \cite[Corollary 2.1]{DDS} for $d=4$, \cite[Theorem 2]{doney1997one} for $d\ge 5$. Equation \eqref{EquationFirst2} can be quickly deduced from the proofs in the aforementioned references, however we give its proof for completeness\footnote{It can alternatively be deduced from the recent preprints \cite{2023arXiv230312505B,2022arXiv220202935V}.}. From the references, we get that for any $\widetilde{\e}>0$ fixed, as $n\to\infty$
\begin{equation}
    \P\rk{S_n=x}=p n x^{-1-d/2}\rk{1+o(1)}=\P\rk{S_n=x, \exists i\in [n]\colon x(1-\widetilde{\e})\le X_i\le x\rk{1+\widetilde{\e}}}\, .
\end{equation}
From this, it follows that, as $n\to\infty$
\begin{equation}
     \P\rk{S_n=x, \nexists i\in[n]\colon x(1-\widetilde{\e})\le X_i\le x\rk{1+\widetilde{\e}}}=o(1) \P\rk{S_n=x}\, .
\end{equation}
It remains to show that the region in between $(1-\widetilde{\e})x$ and $x-a_n^{1+\delta}$ is not likely, same for $\ek{x+a_n^{1+\delta},(1+\widetilde{\e})x}$. We bound using the union bound
\begin{equation}\label{EquationLongBound}
    \P\rk{S_n=x,\exists i\in[n]\colon X_i\in \ek{(1-\widetilde{\e})x,x-a_n^{1+\delta}}}\le n\P\rk{S_{n-1}\ge a_n^{1+\delta}} \sup_{y\in \ek{(1-\widetilde{\e})x,x-a_n^{1+\delta}}}\P\rk{X_1=y}\, .
\end{equation}
Due to Equation \eqref{EquationAssumptionXi}, there exists $C>0$ such that
\begin{equation}
    \sup_{y\in \ek{(1-\widetilde{\e})x,x-a_n^{1+\delta}}}\P\rk{X_1=y}\le C \P\rk{X_1=x}\, .
\end{equation}
Using Equation \eqref{EquationUnifromBound}, we can bound $n\P\rk{X_1=x}$ by $C\P\rk{S_n=x}$, for some (different) $C>0$. Furthermore, by Equation \eqref{EquationUnifromBound}, we can bound $\P\rk{S_n\ge a_n^{1+\delta}}$ by $o(1)$.

Hence, we get that Equation \eqref{EquationLongBound} is bounded from above by
\begin{equation}
    C\P\rk{S_n=x}o(1)
\end{equation}
The same reasoning holds for the upper bound. This concludes the proof.\qed 
\end{proofsect}
We now apply the general results above to our setting: $\rk{\Lcal_x}_x$ takes the role of $\rk{X_i}_i$ and the Poisson property ensures independence. We begin by giving  the precise asymptotics for long loops.
\begin{lemma}\label{LemExpLongLoops}
We have that as $n\to\infty$
\begin{equation}
  M\left[\omega(0)=0,\ell(\omega)=n\right]= \rk{\frac{d}{2\pi }}^{d/2}\frac{1}{n(n\beta)^{d/2}}\oo\, .
\end{equation}
\end{lemma}
\begin{proofsect}{\textbf{Proof of Lemma \ref{LemExpLongLoops}}}
By the definition $M$, we have that
\begin{equation}
    M\left[\omega(0)=0,\ell(\omega)=n\right]=\frac{1}{n}p_{\beta n}(0)=\frac{1}{n}\P_0\rk{\omega(\beta n)=0}\, .
\end{equation}
Recall that by \cite[Theorem 2.1.3]{lawler2010random} as $t\to \infty$
\begin{equation}\label{RWapproximation}
    \P_0(\omega(t)=0)=p_{t}(0)=\rk{\frac{d}{2\pi t}}^{d/2}\rk{1+o(1)}\, .
\end{equation}
This concludes the proof.\qed
\end{proofsect}
We now apply the above results to our setting. For an illustration see Figure \ref{Concentration}. We begin with a corollary to the above results. Recall that $\P_\L$ is the Poisson point process with intensity measure $M_\L$ defined in Equation \eqref{FirstEquationOfM}.
\begin{cor}\label{CorAvLoclExceed}
Recall that $\L=[-N/2,N/2)^d\cap\Z^d$. We have that for $1>c>0$ fixed and all $k$ with $c^{-1}\al\ge k\ge c\al$, as $N\to\infty$
\begin{equation}
    \P_\L\left(\Lcal_\L=k+\rhoc\al\right)\sim\P_\L\left(\exists\omega\colon \omega(0)\in\L\textnormal{ and }\ell(\omega)= k\right)\sim \rk{\frac{d}{2\pi }}^{d/2}\frac{\beta \abs{\L}}{(k\beta)^{1+d/2}}\, ,
\end{equation}
where, as before, the $o(1)$-term implied by the $\sim$ notation does not depend on $k$. 

Furthermore, the random variables $(\Lcal_x)_{x\in\L}$ satisfy a CLT with mean $\rhoc$ and with rate
\begin{equation}\label{equationCLT}
    a_\L=\begin{cases}
    \abs{\L}^{2/3}&\textnormal{ if }d=3\, ,\\
    \sqrt{\abs{\L}\log\abs{\L}}&\textnormal{ if }d=4\, ,\\
    \abs{\L}^{1/2}&\textnormal{ if }d\ge 5\, .\\
    \end{cases}
\end{equation}
\end{cor}
\begin{figure}
    \centering
    \begin{tikzpicture}[very thick]
    \draw[black] (0,0) -- (8,0);
    \draw [decorate,decoration = {calligraphic brace,amplitude=5pt}] (1,0.25) --  (3,0.25)node[pos=0.5,above=5pt,black]{\small{$\Lcal_\L-\ell(\omega)$}};;
    \draw [decorate,decoration = {calligraphic brace, mirror,amplitude=2pt}] (1,-0.65) --  (2,-0.65)node[pos=0.5,below=1pt,black]{\small{$\al^{5/6}$}};
    \draw[thin] (2,-0.075)node[below]{\small{$\rhoc\al$}}  -- (2,0.075);
    \draw[thin] (3,-0.075) -- (3,0.075);
    \draw[thin] (1,-0.075) -- (1,0.075);
    \draw [decorate,decoration = {calligraphic brace,amplitude=5pt}] (5,0.25) --  (7,0.25)node[pos=0.5,above=5pt,black]{\small{$\ell(\omega)$}};;
    \draw[thin] (6,-0.065)node[below]{\small{$\rhoe\al$}}  -- (6,0.065);
    \draw [decorate,decoration = {calligraphic brace, mirror,amplitude=2pt}] (5,-0.75) --  (6,-0.75)node[pos=0.5,below=1pt,black]{\small{$\al^{5/6}$}};
    \draw[thin] (7,-0.075) -- (7,0.075);
    \draw[thin] (5,-0.075) -- (5,0.075);
    \end{tikzpicture}
    \caption{A sketch of the way the total particle number achieves a value of $\rk{\rhoc+\rhoe}\al$: it splits into a single loop $\omega$ around the value $\rho_\mathrm{e}\al$ and the rest concentrates at $\rho_\mathrm{c}\al$, fluctuations being bounded by $\al^{5/6}$.}
    \label{Concentration}
\end{figure}
\begin{proofsect}{\textbf{Proof of Corollary \ref{CorAvLoclExceed}}}
The idea of the proof is to use the fact that conditional on the number of loops, the individual loop lengths are all iid with a nice distribution. We then combine this with the fact that the number of loops follows a Poisson distribution, which concentrates sharply around its mean. Subtracting the mean $\rhoc$ allows us to apply Proposition \ref{PropositionLDP}.

By the fundamental properties of a Poisson point process, we have that
\begin{equation}
    \P_\L\left(\exists\omega\colon \omega(0)\in\L\textnormal{ and }{\ell}(\omega)= k\right)=1-\ex^{-M_\L\left[\ell(\omega)=k\right]}\, .
\end{equation}
Note that by translation invariance
\begin{equation}
     M_\L\left[\ell(\omega)=k\right]=\abs{\L}M_\L\left[\omega(0)=0,\ell(\omega)=k\right]\, .
\end{equation}
Using Lemma \ref{LemExpLongLoops}, we get that, as $N\to\infty$
\begin{equation}
    \P_\L\left(\exists\omega\colon \omega(0)\in\L\textnormal{ and }{\ell}(\omega)= k\right)=1-\ex^{-\abs{\L}M_\L\left[\omega(0)=0,\ell(\omega)=k\right]}\sim\rk{\frac{d}{2\pi }}^{d/2}\frac{\beta\abs{\L}}{(k\beta)^{d/2+1}}\, . 
\end{equation}
We now transfer the result from Proposition \ref{PropositionLDP} to our setting: let $N_x$ be the number of loops starting at $x\in \L$. By the fundamental properties of the Poisson process, we have that $\sum_{x\in\L}N_x$ is distributed like a Poisson random variable with parameter
\begin{equation}
    \abs{\L}\sum_{j\ge 1}\frac{1}{j}p_{\beta j}(0):=\al\Phi(0)\, .
\end{equation}
Furthermore, again by the properties of the Poisson process, conditional on $\sum_{x\in\L}N_x$, the individual loop lengths $\ell$ are iid with law $P$ and distribution
\begin{equation}\label{Equationellk}
    P\rk{\ell=k}=\frac{1}{\Phi(0)}\frac{p_{\beta k}(0)}{k},\quad\text{for}\quad k\ge 1\, .
\end{equation}
By the CLT for iid Poisson random variables, we have that $\sum_{x\in\L}N_x$ will be equal to $\al\Phi(0)$ with Gaussian fluctuations of order $\al^{1/2}$. We enlarge the $1/2$ to $5/6$ to gain some concentration of measure.

Let $E$ be the event the number of loops is typical
\begin{equation}
    E=\gk{\abs{\sum_{x\in\L}N_x-\al\Phi(0)}\le \al^{5/6}}\, .
\end{equation}
We now use a large deviation bound (derived from \cite[Exercise 4.13]{mitzenmacher2017probability}) for Poisson random variables: suppose $X$ is Poisson with parameter $n\l>0$ (where both $n,\l>0$), then there exists $c=c_\l$ such that for all $n>0$
\begin{equation}
    \P\rk{\abs{X-\l n}>\e n}\le 2\ex^{-2nc_\l\e^2}\, .
\end{equation}
Applying the above for $\e=\al^{-1/6}$, we get as $N\to\infty$
\begin{equation}\label{POIssonLDPBound}
        \P\rk{E^c}=\Ocal\rk{\ex^{-\abs{\L}^{1/2}}}\, .
\end{equation}
Hence, the event $E^c$ is negligible. We expand\footnote{here, and from now on, if the sum has non-integer limits, we always take the floor to make the expression well-defined.} using Equation \eqref{Equationellk}
\begin{equation}
    \P_\L\rk{\Lcal_\L=k+\rhoc\al,E}=\sum_{j=-\al^{5/6}}^{\al^{5/6}}\P_\L\rk{\sum_{x\in\L}N_x=\al\Phi(0)+j}P\rk{\sum_{i=1}^{\al\Phi(0)+j}\ell_i=k+\rhoc\al}\, .
\end{equation}
As the $\ell_i$'s satisfy the assumptions from Proposition \ref{PropositionLDP} (see Equation \eqref{Equationellk}), we can now apply it (with $n=\al\Phi(0)+j$) to obtain, as $N\to\infty$
\begin{equation}
    P\rk{\sum_{i=1}^{\al\Phi(0)+j}\!\ell_i=\!k+\rhoc\al}\sim (\al\Phi(0)\!+\!j)P\rk{\ell_1=k}=\frac{\al\Phi(0)+j}{\Phi(0)}\frac{p_{\beta k}(0)}{k}\sim \rk{\frac{d}{2\pi }}^{d/2}\!\!\frac{\beta\al}{(\beta k)^{d/2+1}}\, ,
\end{equation}
using the approximation of the random walk kernel (Equation \eqref{RWapproximation}) in the last step and the fact that $\abs{j}=\Ocal\rk{\al^{5/6}}$. We can plug the above into Equation \eqref{Equationellk} and conclude the proof of the first statement by recalling that $\P_\L\rk{E}=1-o(1)$.

For the second statement, note that tails imply a CLT exactly at the scale stated, as in Lemma \ref{LemmaRefCLT}, see \cite[Section 8.3]{bingham1989regular} for a reference.\qed
\end{proofsect}
The above corollary invites the following intuition, which might be helpful in the following paragraphs: we may ``expand'' (recalling Equation \eqref{MeanLocalTimeEquation})
\begin{equation}
    \sum_{x\in\L}\Lcal_x=\Lcal_\L\approx\rho_\mathrm{c}\abs{\L}+a_\L\mathcal{X}\, ,
\end{equation}
where $\Xcal$ is a stable random variable. This means that the distance between $\Lcalo_\L$ and $\rhoc$ will be of order $a_\L\abs{\L}^{-1}$ with high probability.

Next, we transfer the other results from Proposition \ref{PropositionLDP} to our setting.
\begin{cor}\label{RightCorrrr}
Under the same conditions as in Corollary \ref{CorAvLoclExceed}, we have that as $N\to\infty$
\begin{equation}\label{Equation5161}
    \P_\L\left(\Lcal_\L=k+\rhoc\al\right)\sim\P_\L\left(\Lcal_\L=k+\rhoc\al,\exists \omega\colon \ell(\omega)\in \ek{k-\al^{5/6},k+\al^{5/6}}\right)\, .
\end{equation}
We furthermore have that as $N\to\infty$
\begin{equation}\label{EuqationLargeBOund}
    \P_\L\rk{\Lcal_\L\ge \rhoc\al+\al^{5/6}}=o(1)\, .
\end{equation}
\end{cor}
\begin{proof}
First, note that the CLT scale $a_\L$ of $\rk{\Lcal_x}_{x\in\L}$ grows at most like $\al^{2/3}$. Hence, we can find $\delta>0$ such that $a_\L^{1+\delta}=o\rk{\al^{5/6}}$. We can then employ the LDP-bound in Equation \eqref{POIssonLDPBound} in the same way is in the proof of Corollary \ref{CorAvLoclExceed} to transfer the result (Equation \eqref{EquationFirst2}) from Proposition \ref{PropositionLDP} to the Poissonian setting. This shows that Equation \eqref{Equation5161} holds. The same argument is then used in the proof of Equation \eqref{EuqationLargeBOund}.
\end{proof}
Next we give a bound on the probability that $\Lcal_\L$ will be less than its mean. As the deviations downwards are governed by a standard LDP principle, the decay is quite rapid.
 \begin{proposition}\label{PropositionLowerBound}
For all $N$ large enough
 \begin{equation}
     \P_\L\rk{\Lcal_\L<\rhoc\al-\al^{5/6}}\le {\ex^{-\al^{1/6}}}
 \end{equation}
 \end{proposition}
 \begin{proofsect}{\textbf{Proof of Proposition \ref{PropositionLowerBound}}}
By the Markov inequality, we have that for $t\ge 0$
\begin{equation}\label{Equation33123Markov}
    \P_\L\rk{\Lcal_\L<\rhoc\al-\al^{5/6}}\le \ex^{-t\al^{5/6}}\E\ek{\ex^{t\rk{\Lcal_\L-\rhoc\al}}}=\ex^{-t\al^{5/6}+\al\sum_{j\ge 1}\frac{p_{\beta j}(0)}{j}\rk{\ex^{-tj}-1-tj}}\, ,
\end{equation}
where we use the second statement of Proposition \ref{PropositionPoisson} in the last step. We choose $t=\al^{-1/2}$. We bound, using that as $x\downarrow 0$ we have $\ex^{-x}-1-x=\Ocal\rk{x^2}$,
\begin{equation}
    \sum_{j= 1}^{\al^{1/2}}\frac{p_{\beta j}(0)}{j}\rk{\ex^{-\al^{-1/2}j}-1-\al^{-1/2}j}=\Ocal\rk{\al^{-1}} \sum_{j= 1}^{\al^{1/2}}jp_{\beta j}(0)=\Ocal\rk{\al^{-1+1/4}}\, ,
\end{equation}
using that $d=3$ maximizes the above term. On the other hand
\begin{equation}
    \sum_{j>\al^{1/2}}\frac{p_{\beta j}(0)}{j}\rk{\ex^{-\al^{-1/2}j}-1-\al^{-1/2}j}\le C\al^{-1/2}\sum_{j>\al^{1/2}}{p_{\beta j}(0)}=\Ocal\rk{\al^{-3/4}}\, ,
\end{equation}
using again the fact that $d=3$ maximizes the above. Hence, we get that for $t=\al^{-1/2}$
\begin{equation}
    \al\sum_{j\ge 1}\frac{p_{\beta j}(0)}{j}\rk{\ex^{-tj}-1-tj}=\Ocal\rk{\al^{1/4}}\, .
\end{equation}
Plugging the above into Equation \eqref{Equation33123Markov}, we obtain, as $N\to\infty$
\begin{equation}
    \P_\L\rk{\Lcal_\L<\rhoc\al-\al^{5/6}}\le \ex^{-\al^{5/6-1/2}+\Ocal\rk{\al^{1/4}}}=\ex^{-\al^{1/3}+\Ocal\rk{\al^{1/4}}}\, ,
\end{equation}
and hence the result follows.\qed
 \end{proofsect}
We conclude the section with a small technical lemma on the range of the random walk. Note that $\E_x\ek{\Rcal_t}\sim \kappa_d t$ as $t\to\infty$ (for example \cite{hamana2006}), where $\kappa_d=\P_0\rk{H_{0}=\infty}$ is the escape probability of the simple random walk. We now study the deviations from its mean.
\begin{lemma}\label{LemmaRange}
Let $\Rcal_t$ be the random walk range. Then, there exists $\gamma>0$ such that for every $\e>0$ we have as $t\to \infty$
\begin{equation}\label{EquationRangeEstimate}
    \P_{x,x}^t\rk{\abs{\Rcal_t-\kappa_d t}>\e t}=\Ocal\rk{\ex^{-t^\gamma}}\, ,
\end{equation}
Furthermore, as $t\to\infty$, for $\B_{0,0}^t$ the normalized bridge measure defined in Section \ref{sec_not_setup},
\begin{equation}
    \B_{0,0}^t\ek{\Rcal_t}=\Ocal\rk{t}\qquad\textnormal{and}\qquad\B_{0,0}^t\ek{\Rcal_t^2}=\Ocal\rk{t^2}\, .
\end{equation}
\end{lemma}
\begin{proofsect}{\textbf{Proof of Lemma \ref{LemmaRange}}}
\cite[Theorem 2.5]{hamana2006} and \cite[Theorem 1.2.10]{PP11} give Equation \eqref{EquationRangeEstimate} for the discrete-time random walk. However, the number of jumps for the continuous-time random walk is Poisson distributed. Hence we can use an LDP-bound (see Equation \eqref{POIssonLDPBound}) for Poisson random variables to transfer the result to the continuous-time random walk. 
 
 For the second statement, we can use that $\Rcal_t$ is bounded by the number of jumps up to time $t$ of the random walker. This number is Poisson distributed with parameter $t$ and hence the statement follows.\qed
\end{proofsect}

\subsection{Proof of Theorem \ref{THM1}}
We introduce a new object: recall $\rho=\rhoc+\rhoe$ with $\rhoe>0$. Define
\begin{equation}
    \Re=\gk{j\in\N\colon \rho_\mathrm{e}\abs{\L}-\abs{\L}^{5/6}\le j \le \rho_\mathrm{e}\abs{\L}+\abs{\L}^{5/6}}\, .
\end{equation}
This serves the following purpose: the length of the longest loop varies on the CLT scale from Equation \eqref{equationCLT}. The $\abs{\L}^{5/6}$ dominates $ a_\L$ (for all $d\ge 3$) and this will help to dominate the CLT fluctuations in later parts of the article.

We introduce a new measure $ \overline{M}_{\Delta}$ which gives the distribution of the large loops: for $\Delta=xN+
\L$, $\rho_\mathrm{e}>0$ and $\beta>0$, let $ \overline{M}_{\Delta}$ be defined for an event $A$ as
\begin{equation}\label{Mbardef}
    \overline{M}_{\Delta}[A]=Z^{-1}_\Delta\sum_{x\in\Delta}\sum_{j\in\Re}\frac{1}{j}\P_{x,x}^{\beta j}\rk{A}\P_\L\rk{\Lcal_\L=\rho\al-j}\, ,
\end{equation}
with
\begin{equation}\label{EquationPartFuncDef}
    Z_\Delta=\sum_{x\in\Delta}\sum_{j\in\Re}\frac{1}{j}p_{\beta j}(0)\P_\L\rk{\Lcal_\L=\rho\al-j}\, ,
\end{equation}
The intuition of the definition of $ \overline{M}_{\Delta}$ is as follows: the long loop will, with high probability, have a length $j$ with $j\in\Re$. The remaining loops have a combined particle number of $\rho\al-j$.

Recall that by translation invariance if $\abs{\Delta}=\al$, then the distribution of $\Lcal_\L$ is identical to $\Lcal_\Delta$, and hence also $Z_\L=Z_\Delta$. We hence write $\P_\L\rk{\Lcal_\L=k}$, to keep the notation as concise as possible.

\begin{lemma}\label{LemmaPartitionFunctionFreeCase}
As $N\to\infty$
\begin{equation}\label{EquationPartitionFUnctionFirst}
     Z_\L\sim\rk{\frac{d}{2\pi }}^{d/2}\frac{\beta^{-d/2}}{\rho_\mathrm{e}^{d/2+1}\abs{\L}^{d/2}}\, .
\end{equation}
\end{lemma}
\begin{proof}
First of all, note that uniformly for $j\in\Re$, we have that as $N\to\infty$
\begin{equation}
    p_{\beta j}(0)\sim \rk{\frac{d}{2\pi }}^{d/2}\frac{1}{(\rhoe\al\beta)^{d/2}}\, ,
\end{equation}
compare with Lemma \ref{LemExpLongLoops}. Hence uniformly for $j\in\Re$, we have that as $N\to\infty$
\begin{equation}
    \frac{1}{j}p_{\beta j}(0)\sim \rk{\frac{d}{2\pi }}^{d/2}\frac{\beta^{-d/2}}{(\rhoe\al)^{d/2+1}}\, .
\end{equation}
Note that $\P_{\L}\rk{\Lcal_\L\in \rho\al-\Re}=\P_\L\rk{\abs{\Lcal_\L-\rho_c\al}\le \al^{5/6}}=1-o(1)$ (using Corollary \ref{CorAvLoclExceed}). This concludes the proof.
\end{proof}
Next, we state a series of lemmas which are of technical nature only. They help us to pinpoint which loops contribute asymptotically. As their proofs do not directly link to the random interlacements, we have decided to separate them from the rest, they are contained in Section \ref{SectionApproximationLemmas}.

Recall that $C_N=r_N[-N^{d/2-1}/2,N^{d/2-1}/2)^d\cap\Z^d$ and $\L_N=\bigcup_{x\in C_N}\left(xN+\L\right)$, for $r_N=\log^2(N+1)$.
\begin{lemma}\label{LemAd}
Define $p(x)$ as the weight of the event that there is a loop which starts in the cube centered around $xN$ and intersects $K\subset\Z^d$: 
\begin{equation}
    p(x)=\overline{M}_{xN+\L}\left[\omega\cap K\neq \emptyset\right]\, .
\end{equation}
We then have that, as $\abs{x},N\to\infty$
\begin{equation}
    p(x)=\Ocal\left(\abs{x N}^{2-d}\right)\qquad\textnormal{and}\qquad \sum_{x\in\Z^d}p(x)<\infty\, .
\end{equation}
For the origin $x=0$, we have that $p(0)=\Ocal\left(\abs{N}^{2-d}\right)$.
\end{lemma}
The bound on $p(x)$ is by no means optimal, however it does suffice in our case.

The next lemma shows that the event that a loop outside $\L_N$ contributes to events close to origin, is negligible.
\begin{lemma}\label{LemmaOutsideCN}
For any $K\subset \Z^d$ compact, as $N\to\infty$
\begin{equation}
    \P_N^{\rho}\left(\exists \omega\colon \omega(0)\notin \L_N\,\,\mathrm{ but }\,\,\omega\cap K\neq\emptyset\right)=o(1)\, .
\end{equation}
\end{lemma}

The next lemma proves that for each $x\in C_N$ asymptotically, each box $\Delta=xN+\L$ generates at least one long loop with length in $\Re$:
\begin{lemma}\label{LemmaAtLeastOneLongLoop}
For any $K\subset \Z^d$ compact, as $N\to\infty$
    \begin{equation}
    \P_N^\rho\left(\exists\omega_o\colon \omega_o(0)\in \L_N\textnormal{ and }K\cap\omega_o\neq \emptyset\textnormal{ and }\forall\omega\colon \ell(\omega)\notin\Re\right)=o\left(1\right)\, .
\end{equation}
\end{lemma}

The next lemma proves that the only loop which has a chance to hit a bounded set $K$, is the long loop, unless the loop starts very close to the origin (the case where $\omega(0)\in\L$).
\begin{lemma}\label{LemmaAtMostOneLongLoop}
For any $K\subset \Z^d$ compact, as $N\to\infty$
    \begin{equation}
 \P_N^\rho\left(\exists\omega_o\colon \omega_o(0)\notin\L,\,\omega_o\cap K\neq\emptyset,\,\ell(\omega_o)\notin\Re \right)=o\left(1\right)\, .
\end{equation}
\end{lemma}

The next lemma shows that the probability that a box $xN+\L$ contributes is uniformly approximated by $p(x)$. It furthermore gives that the long loop is distributed according to $\overline{M}_{xN+\L}$. Recall that $\P_{xN+\L}^\rho$ is given by $\P_{xN+\L}$ conditioned on the event $\gk{\Lcal_{xN+\L}=\rho\al}$.

\begin{lemma}\label{LemmaApproximateByPx}
Uniformly in $x\in C_N$, as $N\to\infty$
    \begin{equation}
    \P_{xN+\L}^\rho\rk{\exists \omega\colon\ell(\omega)\in\Re\text{ and }\omega\cap K\neq\emptyset}=p(x)\rk{1+o(1)}\, ,
\end{equation}
as well as
\begin{equation}
    \E_{xN+\L}^\rho\ek{f(\omega)\1\gk{\exists \omega\colon\ell(\omega)\in\Re\text{ and }\omega\cap K\neq\emptyset}}\sim \overline{M}_{xN+\L}[f(\omega),\omega\cap K\neq \emptyset]\, ,
\end{equation}
for $f$ bounded, invariant under path parametrization and local. Uniformly here means that the $o(1)$ term does not depend on $x$. Furthermore, the same asymptotics hold true when replacing $\exists$ by $\exists !$ (there exists a unique).
\end{lemma}
The next lemma shows that for $\Delta=\L$, the process asymptotically normalises:
\begin{lemma}\label{LemmaAtTheOrigin}
    For $F$ a test function as in Definition \ref{definitionLocalconvergence}, we have that as $N\to\infty$
    \begin{equation}
        \E_\L^\rho\ek{F(\eta)}=\E_\L\ek{F\rk{\eta}}\rk{1+o(1)}\, .
    \end{equation}
\end{lemma}


To summarize: we have shown that only loops started in $\L_N$ have a chance to hit a fixed and bound set $K$. Furthermore, any loop hitting $K$ has to either have length in $\Re$ or start very close to the origin. Furthermore, the distribution of the long loops is governed by $\overline{M}_\Delta$. Finally, loops started close to the origin will not be long with high probability.

The next proposition is key, as it allows us to approximate the conditioned process by a free one.
\begin{proposition}\label{LabelLemmaMain}
    Fix $K\subset\Z^d$ bounded and set $\widetilde{\P}_K$ the PPP with intensity measure
\begin{equation}
    \sum_{x\in C_N}\overline{M}_{xN+\L}\left[\,\cdot\,\cap\{\omega\cap K\neq \emptyset\}\right]=:\widetilde{M}_K\, .
\end{equation}
Then for each test function (recall Definition \ref{definitionLocalconvergence}), we have that as $N\to\infty$
\begin{equation}
    \E_N^{\rho}\ek{F\rk{\eta}}\sim \E_{\L}\otimes\widetilde{\E}_K\ek{F}\, .
\end{equation}
\end{proposition}
\begin{proofsect}{\textbf{Proof of Proposition \ref{LabelLemmaMain}}}
The idea is to use the previous lemmas to allow us to reduce $\E_N^\rho$ to just a single long loop in each box $xN+\L$, except for the origin. We then use the representation in Proposition \ref{PropositionPoisson} to prove that the processes is Poisson.
By \cite[Theorem 24.7]{klenke2013probability}, one may assume that
   \begin{equation}\label{Equation13Mar1}
        F(\eta)=\ex^{-\eta[f]}\, ,
   \end{equation}
where $f\colon \Gamma\to [0,\infty]$ and only depends on the values of paths inside the set $K$ and is invariant under reparametrization. 

By Lemma \ref{LemmaOutsideCN} and Lemma \ref{LemmaAtTheOrigin}, we have that loops outside $\L_N$ (or equivalently $xN+\L$ with $x\notin C_N$) do not contribute and that at the origin, the process is not affected by the conditioning. This gives directly
    \begin{equation}
        \E_N^{\rho}\ek{F\rk{\eta}}\stackrel{\text{Lem. \ref{LemmaOutsideCN}}}\sim \rk{\bigotimes_{x\in C_N}\E_{xN+\L}^\rho}\ek{F\rk{\eta}}\stackrel{\text{Lem. \ref{LemmaAtTheOrigin}}}\sim \E_\L\ek{F(\eta)}\rk{\bigotimes_{x\in C_N\setminus\gk{0}}\E_{xN+\L}^\rho}\ek{F\rk{\eta}}\, .
    \end{equation}
    By Lemma \ref{LemmaAtLeastOneLongLoop} and Lemma \ref{LemmaApproximateByPx}, we get that in each box $xN+\L$ there exists precisely one long loop, i.e., 
    \begin{equation}
        \rk{\bigotimes_{x\in C_N\setminus\gk{0}}\E_{xN+\L}^\rho}\ek{F\rk{\eta}}\sim \rk{\bigotimes_{x\in C_N\setminus\gk{0}}\E_{xN+\L}^\rho}\ek{F\rk{\eta}\1\gk{\forall x\in C_N\, \exists !\omega\colon \ell\rk{\omega}\in\Re)}}\, .
    \end{equation}
    Equation \eqref{Equation13Mar1} gives that
    \begin{equation}
        F\rk{\eta\1\gk{\omega(0)\in\L_N}}=\ex^{-\sum_{x\in C_N}\eta[f\1\gk{\omega(0)\in xN+\L}]}\, .
    \end{equation}
    Using Lemma \ref{LemmaApproximateByPx} and the independence of the Poisson processes, we can rewrite the above as
    \begin{multline}
        \prod_{x\in C_N\setminus\gk{0}}\E_{xN+\L}\ek{\ex^{-f(\omega)}\1\gk{\exists! \omega\colon\ell\rk{\omega}\in\Re}}\stackrel{\text{Lem. \ref{LemmaApproximateByPx}}}\sim \prod_{x\in C_N\setminus\gk{0}}\overline{M}_{xN+\L}\ek{\ex^{-f}}\\
        =\prod_{x\in C_N\setminus\gk{0}}\rk{1-p(x)+\overline{M}_{xN+\L}\ek{\ex^{-f}\1\gk{\omega\cap K\neq\emptyset}}}\, .
    \end{multline}
    As for $x=0$ the origin, we have $p(0)=o(1)$, we can rewrite the above as
    \begin{equation}
        \prod_{x\in C_N}\rk{1-p(x)+\overline{M}_{xN+\L}\ek{\ex^{-f}\1\gk{\omega\cap K\neq\emptyset}}}\, .
    \end{equation}
    Denote $q(x)=\overline{M}_{xN+\L}\ek{\ex^{-f}\1\gk{\omega\cap K\neq\emptyset}}$. We expand the product above
    \begin{equation}\label{Equation31420231}
    \sum_{m\ge 0}\sum_{(x_1,\ldots,x_m)\in C_N^m}\rk{\prod_{i=1}^m \frac{q(x_i)}{1-p(x_i)}}\prod_{x\in C_N}\rk{1-p(x)}\, .
    \end{equation}
    The above expansion is exact. Furthermore, as $p(x)$ is summable (see Lemma \ref{LemAd}) and $q(x)\le p(x)$, we have that the above is approximated arbitrarily well by
    \begin{equation}
       \sum_{m=0}^M\sum_{(x_1,\ldots,x_m)\in C_N^m}\rk{\prod_{i=1}^m \frac{q(x_i)}{1-p(x_i)}}\prod_{x\in C_N}\rk{1-p(x)}\, ,
    \end{equation}
    by making $M$ large. Note that by the summability of the $p(x)$, we get that (using Taylor approximation in the second step)
    \begin{equation}\label{Equation3142023}
        \prod_{x\in C_N}\rk{1-p(x)}=\exp\rk{\sum_{x\in C_N}\log\rk{1-p(x)}}\sim \exp\rk{-\sum_{x\in C_N}p(x)}=\exp\rk{-\widetilde{M}_K\ek{\1}}\, .
    \end{equation}
    Note that for $m\le M$, as $N\to\infty$
    \begin{equation}
        \sum_{(x_1,\ldots,x_m)\in C_N^m}\rk{\prod_{i=1}^m \frac{q(x_i)}{1-p(x_i)}}\sim \sum_{(x_1,\ldots,x_m)\in C_N^m}\rk{\prod_{i=1}^m {q(x_i)}}\sim\frac{1}{m!} \sum_{\gk{x_1,\ldots,x_m}\subset C_N}\rk{\prod_{i=1}^m {q(x_i)}}\, ,
    \end{equation}
    as $q(x)\le p(x)$ and $1-p(x)\sim 1$, compare Lemma \ref{LemAd}. By expanding the product, we find that
    \begin{equation}
        \sum_{\gk{x_1,\ldots,x_m}\subset C_N}\rk{\prod_{i=1}^m {q(x_i)}}=\prod_{i=1}^m\rk{\sum_{x\in C_N}q(x_i)}=\prod_{i=1}^m\rk{\widetilde{M}_K\ek{\ex^{-f}}}=\rk{\widetilde{M}_K}^{\otimes m}\ek{\ex^{-f}}\, .
    \end{equation}
    Plugging the last three equations into Equation \eqref{Equation31420231}, we get that, as $N\to\infty$
    \begin{equation}
         \sum_{m\ge 0}\sum_{(x_1,\ldots,x_m)\in C_N^m}\rk{\prod_{i=1}^m \frac{q(x_i)}{1-p(x_i)}}\prod_{x\in C_N}\rk{1-p(x)}\sim \ex^{-\widetilde{M}_K\ek{\1}}\sum_{m\ge 0}\frac{1}{m!}\rk{\widetilde{M}_K}^{\otimes m}\ek{\ex^{-f}}=\widetilde{\E}_K\ek{F}\, ,
    \end{equation}
    where we used the first statement of Proposition \ref{PropositionPoisson} in the last step. This concludes the proof. \qed
\end{proofsect}
Having shown the asymptotic Poisson property, we now set course on identifying the limiting intensity measure. For this, we need to be able to parametrize any elements in $\Gamma^*$ in an unambiguous, deterministic way. We exploit the continuous-time parametrization for this:\\
Given $\omega\in \Gamma_B$ and $K\subset\Z^d$ finite, we can define the parametrization-invariant observable $D_K$, which measures the ``longest time spent at $K$": for a loop $\omega_o\in\Gamma$, let $L_K(\omega_o)$ be the length of the largest (connected) interval $I_K\subset [0,\beta j)$ such that $\omega_o(I_K)\notin K$. We then set
\begin{equation}
    D_K(\omega)=\inf_{\omega_o\in\Pi^{-1}\ek{\Pi(\omega)}}\Big(\beta\ell(\omega)-L_K(\omega_o)\Big)\, ,
\end{equation}
see Figure \ref{FigurePara} for an illustration. As $D_K$ is constant with respect to the shift operator, we can extend it naturally to a random variable on $\Gamma_B^*$.

For $\omega\in W$ a doubly-infinite path, we set $D_K=0$ if $\omega$ does not hit $D_K$. Otherwise, we define $D_K$ to be the duration from the first entrance to the last exit of $K$. $D_K$ does not depend on the parametrization of $\omega$ and thus can be extended to a random variable on $W^*$.

Given $\omega^*\in\Gamma^*$ intersecting a bounded set $K\subset\Z^d$, we define an (almost surely, due to the waiting times being exponentially distributed) unique $\omega\in\Gamma$, denoted by $\amalg_K(\omega^*)$, by requiring
\begin{enumerate}
    \item $\omega(-\infty,0)\notin K$ and $\omega(0)\in K$ if $\omega^*\in W^*$.
    \item $\omega[D_K,\beta j)\notin K$ if $\omega^*\in\Gamma_{B,j}^*$.
\end{enumerate}
See Figure \ref{FigurePara} for an illustration.

\begin{figure}[h]
    \centering
    \begin{tikzpicture}[scale=1.5, very thick]
    \draw[very thick, red](7.6,0.66) -- (8,0.66) node[below]{\tiny{At $K$}};
    \draw[-](0,0) -- (0,0) node[above]{$\omega^*$};
    \draw[densely dotted,thick] (0,0) node[below]{$0$} -- (8,0) node[below]{$\beta j$};
    \draw[-] (0,-0.05) -- (0,0.05);
    \draw[-] (4,-0.05) node[below]{$\beta j/2$} -- (4,0.05);
    \draw[-] (8,-0.05) -- (8,0.05);
    \draw[very thick, red] (5,0) -- (5.3,0);
    \draw[very thick, red] (5.5,0) -- (6,0);
    \draw [decorate,decoration = {calligraphic brace, mirror,amplitude=2pt}] (5,-0.15) --  (6,-0.15)node[pos=0.5,below=2pt,black]{\small{$D_K$}};
    \draw[-](0,-1) -- (0,-1) node[above]{$\amalg_K\rk{\omega^*}$};
    \draw[densely dotted,thick] (0,-1) node[below]{$0$} -- (8,-1) node[below]{$\beta j$};
    \draw[very thick, red] (0,-1) -- (0.3,-1);
    \draw[very thick, red] (0.5,-1) -- (1,-1);
    \draw[-] (1,-0.95) node[below]{$D_K$} -- (1,-1.05);
    \draw[-] (0,-0.95) -- (0,-1.05);
    \draw[-] (4,-0.95) node[below]{$\beta j/2$} -- (4,-1.05);
    \draw[-] (8,-0.95) -- (8,-1.05);
    \draw[-](0,-2) -- (0,-2) node[above]{$\amalg_K\rk{\omega^*}^\mathfrak{s}$};
    \draw[densely dotted,thick] (0,-2) node[below]{$-\beta j/2$} -- (8,-2) node[below]{$\beta j/2$};
    \draw[very thick, red] (4,-2) -- (4.3,-2);
    \draw[very thick, red] (4.5,-2) -- (5,-2);
    \draw[-] (0,-1.95) -- (0,-2.05);
    \draw[-] (4,-1.95) -- (4,-2.05);
    \draw[-] (4,-2)--(4,-2) node[below]{$0$};
    \draw[-] (8,-1.95) -- (8,-2.05);
    \end{tikzpicture}
    \caption{Illustration of the parametrisations: the first one is an arbitrary representative of $\omega^*$. The times at which the loop is at $K$ are thick and red. $\amalg_K(\omega^*)$ is such that the loop is at $K$ in the beginning and then not anymore for the longest possible period of time. $\amalg_K(\omega^*)^\mathfrak{s}$ shifts the parametrisation from $[0,\beta)$ to $[-\beta j/2,\beta j/2)$. Outside the drawn axis, the loop is continued at $\ddagger$.}
    \label{FigurePara}
\end{figure}
Because we parametrize the random interlacements on $(-\infty,\infty)$ and the loops on $[0,\beta j)$, we need to introduce another time shift: let the map $\sfrak\colon [-\beta j/2,\beta j/2)\to \R$ be defined as
\begin{equation}
    \sfrak_j(s)=s+\beta j/2\, .
\end{equation}
For $\omega\in\Gamma_{B,j}$ (i.e., $\ell(\omega)=j$) and $t\in [-\beta j/2,\beta j/2)$, let $\omega^\sfrak(t)=\omega(\sfrak_j(t))$.

We now define a class $\mathfrak{E}$ of events which will be rich enough to capture the behavior of the random interlacements in the topology we chose. As this class has to contain events for both finite loops and infinite interlacements, we need to formulate carefully.
\begin{definition}\label{DefApprx}
A set $E$ belongs to the \textit{approximating class} $\mathfrak{E}$ if it is of the form: $E$ contains all those $\omega^*\in\Gamma^*$ satisfying the following conditions for some compact $K\subset \Z^d$, a sequence of times $-\infty<s_m<\ldots<s_1<0<t_1<\ldots<t_n<\infty$, and a set of points $y_m,\ldots,y_1,x_1,\ldots x_n,z\in\Z^d$
\begin{enumerate}
    \item $\amalg_K(\omega^*)(0)=z$, with $z\in K$,
    \item If $\amalg_K(\omega^*)\in W^*$, then $\amalg_K(\omega^*)({t_i})=x_i$ for $i\in [n]$. If $\amalg_K(\omega^*)\in \Gamma_B^*$, then $\amalg_K(\omega^*)^\sfrak({t_i})=x_i$ for $i\in [n]$.
    \item If $\amalg_K(\omega^*)\in W^*$, then $\amalg_K(\omega^*)({s_i})=y_i$ for $i\in [m]$. If $\amalg_K(\omega^*)\in \Gamma_B^*$, then $\amalg_K(\omega^*)^\sfrak({s_i})=y_i$ for $i\in [m]$.
\end{enumerate}
\end{definition}
In words, an event $E$ belongs to the approximating class if it specifies the finite dimensional distributions relative to a set $K$, after (re-)parameterizing the loops.

We furthermore define the measure $\widetilde{M}_K^*$ on $W^*$ to be $\widetilde{M}_K\circ \Pi$.
The next proposition is key, as it gives the limiting process.
\begin{proposition}\label{LemStTry}
For every $E\in\mathfrak{E}$, we have that
\begin{equation}
    \lim_{N\to\infty}\widetilde{M}_K^*[E]=\rho_\mathrm{e}\nu(E)\, .
\end{equation}
\end{proposition}
\begin{proofsect}{\textbf{Proof of Proposition \ref{LemStTry}}}
We first need to express $\widetilde{M}_K^*[E]$ in terms of $\overline{M}$. Fix points and times such as in Definition \ref{DefApprx}. We rewrite, recalling the definition of $\widetilde{M}_K$ from Proposition \ref{LabelLemmaMain},
\begin{equation}
    \widetilde{M}_K^*[E]= \frac{1}{Z_\L}\sum_{x\in\L_N}\sum_{j\in\Re}\frac{1}{j}\P_{x,x}^{\beta j}\left(H_K<\infty,\, \omega(H_K)=z,\,\omega^\sfrak\circ H_K\in E\right)\P_\L\rk{\Lcal_\L=\rho\al-j}\, .
\end{equation}
Here, we write $\omega^\sfrak\circ H_K$ for $\left(\omega\circ\theta_{H_K}\right)^\sfrak$. In words, we first shift the loop by the time it takes to hit $K$. We then change the domain from $[0,\beta j)$ to $[-\beta j/2,\beta j/2)$. The majority of starting points of the loop are far away from $K$. This causes the loop to locally behave like a random interlacement.

For technical reasons, we have to exclude points close to $K$. However, it is unlikely that the large loop intersecting $K$ starts close to $K$. Define
\begin{equation}
    \L_N^+=\gk{x\in \L_N\colon \abs{x}>\log N}\, .
\end{equation}
\textbf{Claim}: as $N\to\infty$
\begin{equation}\label{EqatuinClaim}
    \frac{1}{Z_\L}\sum_{x\in \L_N\setminus \L_N^+}\sum_{j\in\Re}\frac{1}{j}\P_{x,x}^{\beta j}\left(H_K<\infty,\, \omega(H_K)=z,\,\omega^\sfrak\circ H_K\in E\right)\P_\L\rk{\Lcal_\L=\rho\al-j}=o(1)\, .
\end{equation}
We delay the proof of the claim to the end of the proof.

Recall that $E$ is the event that at the negative times $s_i$, the process is at $y_i$ for $i\in\gk{1,\ldots,m}$, and at the positive times $t_i$, the process is at $x_i$, for $i\in\gk{1,\ldots,n}$. Define $q(E)$, which gives the distribution of the path between the $x_i$'s and $y_i$'s, by
\begin{equation}\label{EqQE}
    q(E)=\left(\prod_{i=0}^{m-1}p_{s_{m-i-1}-s_{m-i}}^K(y_{m-i},y_{m-i-1})\right)\prod_{i=0}^np_{t_{i+1}-t_i}(x_i,x_{i+1})\, ,
\end{equation}
where we apply the convention that $y_0=x_0=z$, $t_0=s_0=0$, and $p_t^K(x,y)=\P_x(\omega(t)=y,H_K>t)$. Denote $T=t_n-s_m$, the total time "spent" at the event $E$. By the strong Markov-property, after time-reversing the random walk, we obtain
\begin{equation}\label{EquationRef}
    \P_{x,x}^{\beta j}\left(H_K<\infty,\, \omega(H_K)=z,\,\omega^\sfrak\circ H_K\in E\right)=q(E)\E_{y_m,x_n}^{\beta j-T}\ek{\1\gk{H_x<\beta j}\1\gk{H_{\L_N^+}<H_K}}\, ,
\end{equation}
see Figure \ref{FigureTime} for an illustration.
\begin{figure}[h]
    \centering
    \begin{tikzpicture}[scale=1.5, thick]
    \draw[very thick, red](7.6,0.66) -- (8,0.66) node[below]{\tiny{At $K$}};
    \draw[-](0,0) -- (0,0) node[above]{$x$};
    \draw[-](8,0) -- (8,0) node[above]{$x$};
    \draw[densely dotted,thick] (0,0) node[below]{$0$} -- (8,0) node[below]{$\beta j$};
    \draw[-] (0,-0.05) -- (0,0.05);
    \draw[-] (8,-0.05) -- (8,0.05);
    \draw[very thick, red] (5,0) -- (5.3,0);
    \draw[very thick, red] (5.5,0) -- (6,0);
    \draw[-] (5,-0.05) -- (5,0.05)node[above]{$z$};
    \draw[-] (6.5,-0.05) -- (6.5,0.05)node[above]{$x_n$};
    \draw[-] (4.5,-0.05) -- (4.5,0.05)node[above]{$y_m$};
    \draw [decorate,decoration = {calligraphic brace, mirror,amplitude=5pt}] (4.5,-0.1) --  (6.5,-0.1)node[pos=0.5,below=3pt,black]{\small{$E$}};
    \draw[-](0,-1) -- (0,-1) node[above]{$y_m$};
    \draw[-](8,-1) -- (8,-1) node[above]{$y_m$};
    \draw[densely dotted,thick] (0,-1) node[below]{\small{$s_m$}} -- (8,-1);
    \draw[-] (0,-0.95) -- (0,-1.05);
    \draw[-] (8,-0.95) -- (8,-1.05);
    \draw[very thick, red] (0.5,-1) -- (0.8,-1);
    \draw[very thick, red] (1,-1) -- (1.5,-1);
    \draw[-] (0.5,-0.95) -- (0.5,0.-1.05);
    \draw[- ](0.5,-1) -- (0.5,-1)node[below]{\small{$0$}};
    \draw[- ](0.5,-1) -- (0.5,-1)node[above]{\small{$z$}};
    \draw[-] (3.5,-0.95) -- (3.5,-1.05);
    \draw[- ](3.5,-1) -- (3.5,-1)node[above]{$x$};
    \draw[-] (2,-0.95)-- (2,-1.05)node[above]{$x_n$};
    \draw[- ](2,-1) -- (2,-1)node[below]{\small{$t_n$}};
    \draw [decorate,decoration = {calligraphic brace, mirror,amplitude=5pt}] (0,-1.3) --  (2,-1.3)node[pos=0.5,below=3pt,black]{\small{$E$}};
    \end{tikzpicture}
    \caption{Illustration of the Equation \eqref{EquationRef}: using the strong Markov property on the bridge from $x$ to $x$ in time $\beta j$, we obtain a bridge from $y_m$ to $x_n$ (time-reversed) during which the walker visits $x$, as well as the finite dimensional distributions as specified by $E$. Times are marked below and locations above the dotted line}
    \label{FigureTime}
\end{figure}
Note that we can expand
\begin{equation}\label{EquationRangeTrick}
\begin{split}
    \sum_{x\in\L_N^+}&\E_{y_m,x_n}^{\beta j-T}\ek{\1\gk{H_x<\beta j}\1\gk{H_{\L_N^+}<H_K}}\\
    &=\E_{y_m,x_n}^{\beta j-T}\ek{\sum_{x\in\L_N^+}\1\gk{H_x<\beta j}\1\gk{H_{\L_N^+}<H_K}}\\
    &=\E_{y_m,x_n}^{\beta j-T}\ek{\#\gk{x\in \L_N^+\colon H_x<\beta j}\1\gk{H_{\L_N^+}<H_K}}\\
    &=p_{\beta j-T}(y_m,x_n)\B_{y_m,x_n}^{\beta j-T}\ek{\#\gk{x\in \L_N^+\colon H_x<\beta j}\1\gk{H_{\L_N^+}<H_K}}
    \end{split}
\end{equation}
We will show that, as $N\to\infty$ (recall $j\in\Re$)
\begin{equation}\label{Equation1632023}
    \B_{y_m,x_n}^{\beta j-T}\ek{\#\gk{x\in \L_N^+\colon H_x<\beta j}\1\gk{H_{\L_N^+}<H_K}}=\kappa_d \beta j \P_{y_m}\rk{H_K=\infty}\rk{1+o(1)}\, ,
\end{equation}
where $\kappa_d=\P_0\rk{H_{0}=\infty}$. The above is of order $j$ times a constant. Hence, error terms of order $o( j)$ are negligible, which we use later.

We estimate the main contribution by first replacing $\#\gk{x\in \L_N^+\colon H_x<\beta j}$ by $\Rcal_{\beta j}$ in Equation \eqref{Equation1632023}, and calculate the error later.\\
By Lemma \ref{LemmaRange}, there exists a $\gamma>0$ such that for every $\e>0$, as $t\to\infty$
\begin{equation}\label{EquationRangeEstimate}
    \P_{x}\rk{\abs{\Rcal_t-\kappa_d t}>\e t}=\Ocal\rk{\ex^{-t^\gamma}}\, ,
\end{equation}
We estimate Equation \eqref{Equation1632023} using Equation \eqref{EquationRangeEstimate}
\begin{equation}
  \E_{y_m,x_n}^{\beta j-T}\ek{\Rcal_{\beta j}\1\gk{H_{\L_N^+}>H_K}}=\kappa_d(\beta  j)p_{\beta j-T}(y_m,x_n)\B_{y_m,x_n}^{\beta j-T}\rk{H_{\L_N^+}>H_K}\rk{1+o(1)}+\Ocal\rk{\ex^{-(\beta j)^\gamma}}\, .
\end{equation}
The stretch-exponential term is negligible, as we sum over $j\in\Re$. As $\dist\rk{H_K,\L_N^+}\to \infty$, we get that
\begin{equation}
    \B_{y_m,x_n}^{\beta j-T}\rk{H_{\L_N^+}>H_K}=\P_{y_m}\rk{H_K=\infty}\rk{1+o(1)}\, .
\end{equation}
Furthermore, note that $p_{\beta j-T}(y_m,x_n)\sim p_{\beta j}(0)$ as $j\to\infty$, as $T$ and $y_m,x_n$ remains fixed. 

We now estimate the error which comes from replacing $\#\gk{x\in \L_N^+\colon H_x<\beta j}$ by $\Rcal_{\beta j}$ in Equation \eqref{Equation1632023}. There are two sources contributing to the error: loops can hit points far away from the origin (outside of $\L_N$) or close to the origin, in (the inside of) $\L_N^+\setminus\L_N$. We bound, using the Cauchy–-Schwarz inequality, the probability that the path exits the domain $\L_N$
\begin{multline}\label{EquationErrorTermm}
    \B_{y_m,x_n}^{\beta j-T}\ek{\#\gk{x\in \L_N^+\colon H_x<\beta j}\1\gk{H_{\L_N^+}>H_K}\1\gk{H_{\L_N^c}<\beta j}}\\
    \le \rk{\B_{y_m,x_n}^{\beta j-T}\ek{\Rcal_{\beta j}^2}}^{1/2}\rk{\B_{y_m,x_n}^{\beta j-T}\rk{H_{\L_N^c}<\beta j}}^{1/2}\, .
\end{multline}
Note that as $j\to\infty$, $\B_{y_m,x_n}^{\beta j-T}\ek{\Rcal_{\beta j}^2}=\Ocal\rk{j^2}$ and $\B_{y_m,x_n}^{\beta j-T}\rk{H_{\L_N^c}<\beta j}=o(1)$, see Lemma \ref{LemmaRange} for the first bound and Equation \ref{EquationReallyGoodBound} for the second. Hence we can bound the above equation by
\begin{equation}
   \Ocal\rk{j^2}^{1/2}o(1)^{1/2}= o\rk{{ j}}\, .
\end{equation}
This shows that Equation \eqref{EquationErrorTermm} is negligible as Equation \eqref{Equation1632023} is of order constant times $j$. Furthermore, the points in $\L_N\setminus\L_N^+$ do not contribute, as $\abs{\L_N\setminus\L_N^+}=\Ocal\rk{\log^d(N)}=o(j)$ for $j\in\Re$. Hence the replacement of $\#\gk{x\in \L_N^+\colon H_x<\beta j}$ by $\Rcal_{\beta j}$ was valid.

To summarize, so far we have shown that, as $N\to\infty$
\begin{equation}
    \sum_{x\in \L_N^+}\P_{x,x}^{\beta j}\left(H_K<\infty,\, \omega(H_K)=z,\,\omega^\sfrak\circ H_K\in E\right)\sim \kappa_d \beta j p_{\beta j}(0)\P_{y_m}\rk{H_K=\infty}q(E)\, .
\end{equation}
As there are no restrictions on $y_m$, we may assume that it is such that $\P_{y_m}\rk{H_K=\infty}\ge 1-\e$, for $\e>0$ arbitrary. As $\P_\L\rk{\Lcal_\L\in \rho\al-\Re}=1-o(1)$, we have that as $N\to\infty$ and $\e\downarrow 0$
\begin{multline}
    \sum_{j\in\Re}\frac{1}{j}\P_\L\rk{\Lcal_\L=\rho\al-j}\sum_{x\in \L_N^+}\P_{x,x}^{\beta j}\left(H_K<\infty,\, \omega(H_K)=z,\,\omega^\sfrak\circ H_K\in E\right)\\
    \sim \kappa_d\beta^{-d/2}q(E)\rhoe^{-d/2}\abs{\L}^{-d/2}\rk{\frac{d}{2\pi}}^{d/2}\, .
\end{multline}
Recall that by Lemma \ref{LemmaPartitionFunctionFreeCase}, as $N\to\infty$
\begin{equation}
    Z_\L\sim\rk{\frac{d}{2\pi }}^{d/2}\frac{\beta^{-d/2}}{\rho_\mathrm{e}^{d/2+1}\abs{\L}^{d/2}}\, .
\end{equation}
Therefore we have that, recalling Equation \eqref{EquationDefInterl}, as $N\to\infty$
\begin{equation}
    Z^{-1}_\L\sum_{x\in\L_N^+}\sum_{j\in\Re}\P_\L\rk{\Lcal_\L=\rho\al-j}\frac{1}{j}\P_{x,x}^{\beta j} \left(H_K<\infty,\, \omega(H_K)=z,\,\omega^\sfrak\circ H_K\in E\right)\sim \rho_\mathrm{e}\mathtt{Q}^K(E)\, .
\end{equation}
This concludes the proof, modulo the claim in Equation \eqref{EqatuinClaim}. However, note that
\begin{multline}
    \sum_{x\in \L_N\setminus \L_N^+}\sum_{j\in\Re}\P_\L\rk{\Lcal_\L=\rho\al-j}\frac{1}{j}\P_{x,x}^{\beta j}\left(H_K<\infty,\, \omega(H_K)=z,\,\omega^\sfrak\circ H_K\in E\right)\\
    =\Ocal\rk{\log^d(N)\al^{-d/2-1}}\, ,
\end{multline}
and the claim follows upon recalling that $Z_\L=\Ocal\rk{\al^{d/2}}$.
\end{proofsect}
The next lemma is a consequence of the Proposition \ref{LemStTry}.
\begin{lemma}\label{LemTrueConv}
Suppose that $F$ is a path-local, bounded, and continuous function on $\Gamma$ which is also shift-invariant. We then have that
\begin{equation}
    \lim_{N\to\infty}\widetilde{M}_K^*[F]=\rho_\mathrm{e} \nu(F)\, .
\end{equation}
\end{lemma}
\begin{proofsect}{\textbf{Proof of Lemma \ref{LemTrueConv}}}
Due to the properties prescribed onto $F$, it suffices to show that $\widetilde{M}_K^*\circ \amalg_K^{-1}$ converges weakly to $\nu\circ \amalg_K^{-1}$ on $D(-\infty,\infty)$, the space of doubly infinite cadlag paths with values in $\Z^d$. However, this follows easily: using \cite[Theorem 16.8]{billingsley1968convergence}, it suffices to show the convergence of the finite dimensional distributions together with a tightness criterion. The convergence of the finite dimensional distributions was shown in Proposition \ref{LemStTry}. The tightness criterion is equivalent to showing that
\begin{equation}\label{Equation2161}
    \lim_{a\to\infty}\limsup_{N\to\infty}\widetilde{M}_K^*\circ \amalg_K^{-1}\left[\exists t\in [-M,M]\colon \abs{\omega(t)}>a\right]=0\, ,
\end{equation}
for every $M>0$. However we can bound using the same reasoning as in the proof of Proposition \ref{LemStTry}, Equation \eqref{Equation2161} from above by
\begin{equation}
\begin{split}
    \widetilde{M}_K^*\circ &\amalg_K^{-1}\left[\exists t\in [-M,M]\colon \abs{\omega(t)}>a\right]\\
    &\le C\abs{\L}^{d/2}\sum_{z\in K}\sum_{j\in\Re}\P_{z,z}^{\beta j}\left(\exists  t\in [0,M]\cup[\beta j-M,\beta j]\textnormal{ with }\abs{\omega(t)}>a\right)\, .
    \end{split}
\end{equation}
Using Equation \eqref{EquationReallyGoodBound}, the above can be bounded by $\Ocal\left(\ex^{-c a/M}\right)$ uniformly in $N$, for some $c>0$. This concludes the proof.\qed 
\end{proofsect}
\begin{proofsect}{\textbf{Proof of Theorem \ref{THM1}} }
Theorem \ref{THM1} now follows from Lemma \ref{LemTrueConv} and Proposition \ref{LabelLemmaMain} together with the previous results.
Fix $F$ satisfying the conditions of Theorem \ref{THM1} and be supported on $K$. We then have by Proposition \ref{LabelLemmaMain} that for $\rho_\mathrm{e}>0$
\begin{equation}
    \E^{\rho}_N[F]=\E_\L\otimes\widetilde{\E}_K[F]\left(1+o(1)\right)\, .
\end{equation}
where $\widetilde{\E}_K[F]$ is the expectation with respect to $\widetilde{\P}_K[F]$, defined in Proposition \ref{LabelLemmaMain}. However, $\widetilde{\P}_K$ is a PPP with uniformly (in $N$) bounded intensity measure $\widetilde{M}_K$. Since $\widetilde{M}_K$ has uniformly bounded total mass, we have that the convergence of the intensity measure implies the convergence of the Poisson point process, see Proposition \ref{PropositionPoisson}. Since $\P_\L=\P_{\L,\beta,0}\to\P_{\Z^d,\beta,0}$ as $N\to\infty$, we get that
\begin{equation}
    \lim_{N\to\infty}\E^{\rho}_N[F]=\lim_{N\to\infty}\E_{\L,\beta,0}\otimes\widetilde{\E}_K[F]=\E_{\Z^d,\beta,0}\otimes{\E}^{\iota}_{\rhoe}[F]\, .
\end{equation}
This concludes the proof.\qed 
\end{proofsect}
\subsection{Proof of the technical lemmas}\label{SectionApproximationLemmas}
In this section, we prove the various technical lemmas from the previous section.
\begin{proofsect}{\textbf{Proof of Lemma \ref{LemAd}}}
Our goal is to show for bounded $K\subset\Z^d$, as $N\to\infty$
\begin{equation}
    p(x)=\Ocal\left(\abs{x N}^{2-d}\right)\qquad\textnormal{and}\qquad \sum_{x\in\Z^d}p(x)<\infty\, ,
\end{equation}
where
\begin{equation}
    p(x)=\overline{M}_{xN+\L}\left[\omega\cap K\neq \emptyset\right]\, .
\end{equation}
Note that for some $C>0$ (see \cite[Lemma 2.2]{QPhD20} for a direct statement or \cite[Proposition 2.4.5]{lawler2010random} for the underlying large deviation estimate)
\begin{equation}\label{EquationBundFarAway}
    \P_{y,y}^{\beta j}\left(H_K<\beta j\right)=\Ocal\left(j^{-d/2}\abs{y}^{2-d}\Gamma\left(d/2-1,C\abs{y}^2/j\right)\right)\, ,
\end{equation}
as $j,\abs{y}\to\infty$. Here, we use the notation $\Gamma(s,a)$ for the incomplete gamma function, defined as $\Gamma(s,a)=\int_a^\infty t^{s-1}\ex^{-t}\d t$. We employ the rough bound $\Gamma(s,a)=\Ocal(1)$, for $s>0$. Thus $\P_{y,y}^{\beta j}\left(H_K<\beta j\right)=\Ocal\left(j^{-d/2}\abs{y}^{2-d}\right)$ in Equation \eqref{EquationBundFarAway}, uniformly in $j,y$.

Recall that
 \begin{equation}
     p(x)=\overline{M}_{xN+\L}\left[\omega\cap K\neq \emptyset\right]=\frac{1}{Z_\L}\sum_{y\in xN+\L}\sum_{j\in\Re}\frac{1}{j} \P_{y,y}^{\beta j}\left(H_K<\beta j\right)\P_\L\rk{\Lcal_\L=\rho\al-j}\, .
 \end{equation}
Using the above, we can bound this from above by
\begin{equation}
    \frac{1}{Z_\L}\sum_{y\in xN+\L}\sum_{j\in\Re}\Ocal\left(j^{-d/2-1}\abs{y}^{2-d}\right)\P_\L\rk{\Lcal_\L=\rho\al-j}\, .
\end{equation}
Observe that, as $N\to\infty$
\begin{equation}
    \sup_{j\in\Re} \frac{1}{j}\Ocal\left(j^{-d/2}\right)=\Ocal\left(\abs{\L}^{-d/2-1}\right)\quad\textnormal{and}\quad\sum_{y\in xN+\L}\Ocal\left(\abs{y}^{2-d}\right)=\Ocal\left(\abs{x}^{2-d}N^2\right)\, .
\end{equation}
Recall that as $N\to\infty$, we have $\P_\L\rk{\Lcal_\L\in \rho\al-\Re}\sim 1$. By Equation \eqref{EquationPartitionFUnctionFirst}, we have that $Z_\L^{-1}=\Ocal\rk{\abs{\L}^{d/2}}$, and hence as $N\to\infty$
\begin{equation}
    p(x)= \abs{\L}^{d/2}\Ocal\left(\abs{\L}^{-d/2-1}\abs{x}^{2-d}N^2\right)=\Ocal\left(\abs{x N}^{2-d}\right)\, ,
\end{equation}
using $\abs{\L}=N^d$. This concludes the proof of the first part.

For the second part, recall the bound in Equation \eqref{EquationBundFarAway}. Note that also that for each $s>0$, there exists $C>0$ such that for all $x>0$
    \begin{equation}
       \Gamma\rk{s,x}= \int_x^\infty t^{s-1}\ex^{-t}\d t\le C \ex^{-x/2}\, .
    \end{equation}
    Hence, for some different $C>0$, we can bound
    \begin{equation}\label{EquationReallyGoodBound}
        \frac{1}{j}\P_{y,y}^{\beta j}\rk{H_K<\beta j}=\Ocal\rk{{j^{-d/2-1}\abs{y}^{2-d}\ex^{-C\abs{y}^2/j}}}\, .
    \end{equation} 
We bound $\rhoe\al/2\le j\le 2\rhoe\al$ and hence get $j^{-d/2-1}=\Ocal\rk{\al^{-1-d/2}}$. We obtain with the previous equation
\begin{equation}
    \begin{split}
        \sum_{x\in\Z^d}p(x)&=\Ocal\rk{\al^{d/2}}\sum_{y\in\Z^d}\Ocal\rk{\al^{-d/2-1}}\abs{y}^{2-d}\ex^{-C\abs{y}^2/\al}\\
        &\le \Ocal\rk{\al^{-1}}\int_0^\infty r^{2-d+d-1}\ex^{-Cr^2/\al}\d r\\
        &\le \Ocal\rk{1}\int_0^\infty r\ex^{-Cr^2}\d r<\infty\, ,
    \end{split}
\end{equation}
where we bounded the sum by an integral and then employed a change of variables. This concludes the proof.\qed
\end{proofsect}
\begin{proofsect}{\textbf{Proof of Lemma \ref{LemmaOutsideCN}}}
Our goal is to show that for any $K\subset \Z^d$ bounded
\begin{equation}
    \P_N^{\rho}\left(\exists \omega\colon \omega(0)\notin \L_N\,\,\mathrm{ but }\,\,\omega\cap K\neq\emptyset\right)=o(1)\, .
\end{equation}
    We prove this by employing a union bound and then a case distinction with regards to the length of the loop $\omega$. Crucial again are the estimates on long loop lengths from Corollary \ref{CorAvLoclExceed}

    Recall the definition of $\P_N^\rho$ from Equation \eqref{EquationDefinitionrhon}. By a union bound, we get that
    \begin{equation}\label{Equation1471}
        \P_N^{\rho}\left(\exists \omega\colon \omega(0)\notin \L_N\,\,\mathrm{ but }\,\,\omega\cap K\neq\emptyset\right)\le \sum_{x\notin C_N}\P_{xN+\L}\rk{\exists \omega\colon \omega\cap K\neq\emptyset|\Lcal_{xN+\L}=\rho\al}
    \end{equation}
Next, we do expand the above by distinguishing the length of $\omega$: for $\Delta=xN+\L$, using the fourth statement from Proposition \ref{PropositionPoisson}, we get that
    \begin{equation}\label{Bound134}
        \P_\Delta\rk{\exists \omega\colon \omega\cap K\neq\emptyset,\Lcalo_\Delta=\rho}\le \sum_{j=1}^{\rho\al}M_\Delta\ek{\ell(\omega)=j,\,\omega\cap K\neq\emptyset}\P_\Delta\rk{\Lcal_\Delta=\rho\al-j}\, .
    \end{equation}
Note that $\P_\Delta\rk{\Lcal_\Delta=k}=\P_\L\rk{\Lcal_\L=k}$ and hence the probability of the conditioned event on the right-hand side of Equation \eqref{Equation1471} does not depend on $x$.

Our aim is to show that summing the above over $x\notin C_N$, the resulting term is $o\rk{\al^{-d/2}}$. We split the sum in $j\le \delta\abs{\L}$ with $0<\delta<\rhoe$, and $j>\delta\abs{\L}$.
    
\textbf{The case }$j\le \delta\abs{\L}$: in this case, we use that the probability of the particle number being $\rho\al-j$ is bounded by $\Ocal\rk{\al^{-d/2}}$, as $\rho\al-j$ is bounded away from the mean $\rhoc\al$. We estimate
    \begin{equation}\label{Equation3929231}
        \sum_{j=1}^{\delta\abs{\L}}M_\Delta\ek{\ell(\omega)=j,\,\omega\cap K\neq\emptyset}\P_\Delta\rk{\Lcal_\Delta=\rho\al-j}\le C\abs{\L}^{-d/2}M_\Delta\ek{\ell(\omega)\le\delta\abs{\L},\,\omega\cap K\neq\emptyset}\, ,
    \end{equation}
    as for $j$ as above, we have $\rho\al-j\ge \rk{\rhoc+\rhoe-\delta}\al$ and hence $\P_\Delta\rk{\Lcal_\Delta=\rho\al-j}=\Ocal\rk{\abs{\L}^{-d/2}}$ by Corollary \ref{CorAvLoclExceed}. Recall that
    \begin{equation}
        M_\Delta\ek{\ell(\omega)\le\delta\abs{\L},\,\omega\cap K\neq\emptyset}=\sum_{y\in xN+\L}\sum_{j=1}^{\delta\abs{\L}}\frac{1}{j}\P_{y,y}^{\beta j}\rk{H_K<\beta j}\, ,
    \end{equation}
and that $\cup_{x\in C_N}\rk{xN+\L}=\L_N$. We can estimate the above using Equation \eqref{EquationReallyGoodBound}
    \begin{equation}
        \sum_{y\notin\L_N}\sum_{j=1}^{\delta\abs{\L}} \frac{1}{j}\P_{y,y}^{\beta j}\rk{H_K<\beta j}\le  \sum_{y\notin\L_N}\sum_{j=1}^{\delta\abs{\L}} \frac{1}{j}\Ocal\rk{j^{-d/2}\abs{y}^{2-d}\ex^{-C\abs{y}^2/j}}\le C \sum_{y\notin\L_N}\abs{y}^{2-d}\ex^{-C\abs{y}^2/(\delta\abs{\L})}\, .
    \end{equation}
Bounding the sum by an integral, we obtain after a change of variables
\begin{equation}
    \sum_{y\notin\L_N}\abs{y}^{2-d}\ex^{-C\abs{y}^2/(\delta\abs{\L})}\le C\int_{r_N\al}^\infty t^{2-d+d-1}\ex^{-C\abs{t}^2/(\delta\abs{\L})}\d t\le C \abs{\L}\int_{r_N}^\infty t \ex^{-t^2}\d t=o(1)\, .
\end{equation}
Combining the above with Equation \eqref{Equation3929231}, we obtain that
\begin{equation}
      \sum_{y\notin\L_N}\sum_{j=1}^{\delta\abs{\L}}M_\Delta\ek{\ell(\omega)=j,\,\omega\cap K\neq\emptyset}\P_\Delta\rk{\Lcal_\Delta=\rho\al-j}\le o(1)\Ocal\rk{\al^{-d/2}}= o(1)\P_\Delta\rk{\Lcalo_\Delta=\rho}\, .
\end{equation}
Hence,
\begin{equation}
    \sum_{x\notin C_N}\P_{xN+\L}\rk{\exists \omega\colon \omega\cap K\neq\emptyset\text{ and }\ell(\omega)\le \delta\al|\Lcal_{xN+\L}=\rho\al}=o(1)\, .
\end{equation}
\textbf{The case} $j\in\ek{\delta\al,\rho\al}$: for the case $j\in\ek{\delta\al,\rho\al}$, note that using Equation \eqref{EquationReallyGoodBound}, we have that, approximating the sum by and integral
\begin{multline}
    \sum_{y\notin\L_N}M_\Delta\ek{\ell(\omega)=j,\,\omega\cap K\neq\emptyset}= \sum_{y\notin\L_N}\Ocal\rk{j^{-d/2-1}\abs{y}^{2-d}\ex^{-C\abs{y}^2/j}}\\
    =\Ocal\rk{j^{-d/2-1}}\int_{r_N\al}^\infty t\ex^{-C\abs{t}^2/(\delta\abs{\L})}\d t\, ,
\end{multline}
as $N\to\infty$. Changing variables ($t\mapsto t\sqrt{\al}$) and recalling that $j=\Ocal\rk{\al}$, we obtain
\begin{equation}
    \sup_{j\in\ek{\delta\abs{\L},\rho\al}}  \sum_{y\notin\L_N}M_\Delta\ek{\ell(\omega)=j,\,\omega\cap K\neq\emptyset}=\Ocal\rk{\abs{\L}^{-d/2}\ex^{-Cr_N^2}}\, .
\end{equation}
 Recall that
\begin{equation}
    \P_\Delta\rk{\Lcalo_\Delta\in [0,\rhoc+\rhoe-\delta]}=1+o(1)\, .
\end{equation}
Hence, plugging the two previous equations into Equation \eqref{Bound134}, we obtain
\begin{multline}
     \sum_{y\notin\L_N} \sum_{j=\delta\abs{\L}}^{\rho\abs{\L}}\P_\Delta\rk{\Lcalo_\Delta=\rho\al-j}M_\Delta\ek{\ell(\omega)=j,\,\omega\cap K\neq\emptyset}=\Ocal\rk{\abs{\L}^{-d/2}\ex^{-Cr_N^2}}     =o\rk{\abs{\L}^{-d/2}}\, .
\end{multline}
This concludes the proof of the lemma.\qed
\end{proofsect}

\begin{proofsect}{\textbf{Proof of Lemma \ref{LemmaAtLeastOneLongLoop}}}
Our goal is to show
\begin{equation}\label{EquationStatementLemma412}
    \P_N^\rho\left(\exists\omega_o\colon \omega_o(0)\in \L_N\text{ and }K\cap\omega_o\neq \emptyset\text{ and }\forall\omega\colon \ell(\omega)\notin\Re\right)=o\left(1\right)\, .
\end{equation}
    To prove the lemma, a fine control on the probability of long loops is needed.
    We begin by employing a union bound
    \begin{equation}
        \eqref{EquationStatementLemma412}\le  \sum_{x\in C_N}\P_{xN+\L}\left(\exists\omega_o\colon K\cap\omega_o\neq \emptyset\text{ and }\forall\omega\colon \ell(\omega)\notin\Re|\Lcalo_{xN+\L}=\rho\right)
    \end{equation}
For $x\in C_N$, we write $\Delta=xN+\L$. Using the same method as in the estimation of Equation \eqref{EquationRangeTrick}, we can bound
\begin{equation}\label{EquationEstimateByRange}
    \begin{split}
        \sum_{x\in C_N}M_\Delta[\omega\colon\ell(\omega)=k,\, 0\in\omega]=&\sum_{x\in C_N}\sum_{y\in\Delta}\frac{1}{k}\E_{y,y}^{\beta k}\left[\1\{H_0<\beta j\}\right]\\
        \le &\sum_{y\in\L_N}\frac{1}{k}\E_{0,0}^{\beta j}\left[\1\{H_y<\beta k\}\right]\\
        \le &\frac{1}{k}\E_{0,0}^{\beta j}\left[\Rcal_{\beta k}\right]\le C k^{-1}p_{k\beta}(0)\B_{0,0}^{\beta j}\left[\Rcal_{\beta k}\right]\le C k^{-d/2}\, ,
    \end{split}
\end{equation}
using Lemma \ref{LemmaRange} in the last step. Abbreviate $\gk{\forall\omega\colon \ell(\omega)\notin\Re}=\Oe^c$. We expand using fourth statement of Proposition \ref{PropositionPoisson} as in Equation \eqref{Bound134}
\begin{equation}
    \begin{split}
        \P_\Delta\rk{\exists\omega_o\colon K\cap\omega_o\neq \emptyset,\Oe^c,\Lcalo=\rho}&\le \sum_{k\in\N\setminus \Re}\P_\Delta\rk{\exists \omega\colon \omega\cap K\neq \emptyset,\, \ell(\omega)=k,\,\Oe^c,\Lcalo=\rho}\\
        &\le\sum_{k\in\N\setminus \Re}M_\Delta\ek{\omega\cap K\neq \emptyset,\, \ell(\omega)=k}\P_\Delta\rk{\Oe^c,\Lcal_\Delta=\rho\abs{\L}-k}\, .
    \end{split}
\end{equation}
Combining the above with Equation \ref{EquationEstimateByRange}, we obtain the bound
\begin{equation}
    \sum_{x\in C_N} \P_\Delta\rk{\exists\omega_o\colon K\cap\omega_o\neq \emptyset,\Oe^c,\Lcalo=\rho}\le C\sum_{k\in\N\setminus \Re}k^{-d/2}\P_\Delta\rk{\Oe^c,\Lcal_\Delta=\rho\abs{\L}-k}\, .
\end{equation}
We employ different estimates based on the value of $k$.

\textbf{The case} $k\le \al^{\delta}$: we first treat the case very $k$ small. Note that by Corollary \ref{RightCorrrr}, we have that for some $\delta>0$
\begin{equation}
    \sum_{k\le \al^{\delta}}k^{-d/2}\P_\Delta\rk{\Oe^c,\Lcal_\Delta=\rho\abs{\L}-k}=\sum_{k\le \al^{\delta}}k^{-d/2}o(1)\P_\Delta\rk{\Lcal_\Delta=\rho\abs{\L}-k}=o\rk{\al^{-d/2}}\, ,
\end{equation}
as $N\to\infty$. This concludes the proof for the case $k\le \al^{\delta}$.

\textbf{The case} $\al^{\delta}\le k\le \rhoe\al/2$: on the other hand, we have that for $\al^\delta\le k\le \rhoe\al/2$ that $\P_\Delta\rk{\Oe^c,\Lcal_\Delta=\rho\abs{\L}-k}\le \P_\Delta\rk{\Lcal_\Delta=\rho\abs{\L}-k}$, where the latter is bounded by $\Ocal\rk{\al^{-d/2}}$, using Corollary \ref{CorAvLoclExceed}. Hence,
\begin{equation}
    \sum_{k=\al^{\delta}}^{\rhoe\al/2}k^{-d/2}\P_\Delta\rk{\Oe^c,\Lcal_\Delta=\rho\abs{\L}-k}\le \Ocal\rk{\al^{-d/2}}\sum_{k=\al^{\delta}}^{\rhoe\al/2}k^{-d/2}=o\rk{\al^{-d/2}}\, .
\end{equation}
\textbf{The case} $k\ge \rhoe\al/2$: finally, we have using Corollary \ref{CorAvLoclExceed} and Proposition \ref{PropositionLowerBound} and using a supremum bound on $k^{-d/2}$
\begin{equation}
    \sum_{\mycom{k\in\N\setminus \Re}{k\ge \rhoe\abs{\L}/2}}k^{-d/2}\P_\Delta\rk{\Oe^c,\Lcal_\Delta=\rho\abs{\L}-k}\le \Ocal\rk{\al^{-d/2}}\P_\Delta\rk{{\Lcal_\Delta-\rhoc\al}>\al^{5/6}}=o\rk{\al^{-d/2}}\, .
\end{equation}
This concludes the proof.\qed
\end{proofsect}
\begin{proofsect}{\textbf{Proof of Lemma \ref{LemmaAtMostOneLongLoop}}}
Using the union bound as in the lemmas before, we find that statement follows from
\begin{equation}
    \sum_{x\in C_N\setminus\gk{0}} \P_\Delta^\rho\left(\exists\omega_o\colon \omega_o\cap K\neq\emptyset,\,\ell(\omega_o)\notin\Re \right)=o\left(1\right)\, ,
\end{equation}
We show the above by first expanding for $x\in C_N\setminus\gk{0}$ and $\Delta=xN+\L$
    \begin{equation}\label{EquationEq4181}
        \P_\Delta\left(\exists\omega_o\colon K\cap\omega_o\neq\emptyset,\,\ell(\omega_o)\notin \Re,\,\Lcalo_\Delta=\rho\right)=\sum_{j\in \N\setminus\Re}\P_\Delta\left(\exists\omega_o\colon K\cap\omega_o\neq\emptyset,\,\ell(\omega_o)=j,\,\Lcalo_\Delta=\rho\right)\, .
    \end{equation}
    By the last statement of Proposition \ref{PropositionPoisson}, we have that
    \begin{equation}
        \P_\Delta\left(\exists\omega_o\colon K\cap\omega_o\neq\emptyset,\,\ell(\omega_o)=j,\,\Lcalo_\Delta=\rho\right)\le M_{\Delta}[\ell(\omega)=j,\omega\cap K\neq\emptyset]\P_\Delta\left(\Lcal_\Delta=\rho\abs{\L}-j\right)\, .
    \end{equation}
    By Equation \eqref{EquationEstimateByRange}, we have that 
    \begin{equation}
        \sum_{x\in C_N}M_{xN+\L}[\ell(\omega)=j,\omega\cap K\neq\emptyset]= \sum_{x\in C_N}M_{\Delta}[\ell(\omega)=j,\omega\cap K\neq\emptyset]=\Ocal\rk{j^{-d/2}}\, .
    \end{equation}
    We first treat the case $j$ smaller than $\rhoe\abs{\L}-\al^{5/6}$ in Equation \eqref{EquationEq4181}.
    
    \textbf{Case 1}: estimating the case where $j\le \rhoe\al-\al^{5/6}$. We estimate (using Corollary \ref{RightCorrrr} in the last step)
    \begin{equation}
        \sum_{j=\rhoe\abs{\L}/2}^{\rhoe\abs{\L}-\al^{5/6}}j^{-d/2}\P_\Delta\left(\Lcal_\Delta=\rho\abs{\L}-j\right)\le \Ocal\rk{\abs{\L}^{-d/2}}\P_\Delta\left(\Lcal_\Delta>\rhoc\abs{\L}+\abs{\L}^{5/6}\right)=o\rk{\abs{\L}^{-d/2}}\, .
    \end{equation}
    On the other hand, using Corollary \ref{CorAvLoclExceed}
    \begin{equation}
        \sum_{j=\abs{\L}^{5/6}}^{\rhoe\abs{\L}/2}j^{-d/2}\P_\Delta\left(\Lcal_\Delta=\rho\abs{\L}-j\right)\le \sup_{k\ge \rk{\rhoc+\rhoe/2}\al}\P\rk{\Lcal_\Delta=k}\sum_{j=\abs{\L}^{5/6}}^{\rhoe\abs{\L}/2}j^{-d/2}=o\rk{\abs{\L}^{-d/2}}\, ,
    \end{equation}
    since for $k$ as above we have $\P\rk{\Lcal_\Delta=k}=\Ocal\rk{\al^{-d/2}}$ and $j^{-d/2}$ is summable.
    
    Finally, note that by Equation \eqref{EquationReallyGoodBound} as in the proof of Lemma \ref{LemmaOutsideCN}
    \begin{equation}
         \sum_{j=1}^{\abs{\L}^{5/6}}\sum_{x\in C_N\setminus\gk{0}}\P_{x,x}^{\beta j}\rk{H_K<\beta j}=o(1)\, .
    \end{equation}
    This gives using Proposition \ref{PropositionLDP}
    \begin{multline}
        \sum_{j=1}^{\abs{\L}^{5/6}}\sum_{x\in C_N\setminus\gk{0}}M_{xN+\L}[\ell(\omega)=j,\omega\cap K\neq\emptyset]\P_\Delta\left(\Lcal_\Delta=\rho\abs{\L}-j\right)\\\le o(1)\sup_{k\ge 0}\P\rk{\Lcal_\Delta=\rho\al -\abs{\L}^{5/6}+k}=o\rk{\abs{\L}^{-d/2}}\, .
    \end{multline}
     \textbf{Case 2}: we now treat the case $j$ larger than $\rhoe\abs{\L}+\al^{5/6}$. Using Proposition \ref{PropositionLowerBound}, we obtain that
     \begin{equation}
          \P_\Delta\left(\exists\omega_o\colon K\cap\omega_o\neq\emptyset,\,\ell(\omega_o)\ge \rhoe\al+\al^{5/6},\,\Lcalo_\Delta=\rho\right)\le \ex^{-\al^{1/6}}\, .
     \end{equation}
     As the exponential decay dominates the polynomial size of $C_N$, we get that
     \begin{equation}
         \sum_{x\in C_N\setminus\gk{0}}\P_\Delta\left(\exists\omega_o\colon K\cap\omega_o\neq\emptyset,\,\ell(\omega_o)\ge \rhoe\al+\al^{5/6},\,\Lcalo_\Delta=\rho\right)=o\rk{\abs{\L}^{-d/2}}\, .
     \end{equation}
    This concludes the proof.\qed
\end{proofsect}
\begin{proofsect}{\textbf{Proof of Lemma \ref{LemmaApproximateByPx}}}
First note that the probability that there exists two or more loops with length in $\Re$ is negligible: if there were two or more loops with length in $\Re$, the remaining loops would have to have a length of at most $\rhoc\al-\rhoe\al+2\al^{5/6}$. However, by Proposition \ref{PropositionLowerBound}, this even has stretch-exponentially small weight. It is hence negligible.

This allows us to write
\begin{equation}
    \P_N^\rho\rk{\forall x\in C_N\,\,\exists \omega\colon \ell\rk{\omega}\in\Re}\sim \P_N^\rho\rk{\forall x\in C_N\,\,\exists! \omega\colon \ell\rk{\omega}\in\Re}\, .
\end{equation}
Note that by the fourth statement of Proposition \ref{PropositionPoisson}, we have that
\begin{multline}
    \P_{xN+\L}\rk{\exists!\omega\colon \ell(\omega)\in\Re,\, \omega\cap K\neq\emptyset,\, \Lcalo=\rho}\\=\sum_{j\in\Re}M_{xN+\L}\ek{\ell(\omega)=j, \omega\cap K\neq\emptyset} \P_{xN+\L}\rk{\Lcal=\rho\al -j}\, .
\end{multline}
Comparing this with Equation \eqref{Mbardef} concludes the first statement. The second statement follows in the same way.
\end{proofsect}

\begin{proofsect}{\textbf{Proof of Lemma \ref{LemmaAtTheOrigin}}}
This lemma follows from \cite{armendariz2011conditional}, which characterises the distribution of the short loops: to realise the large deviation event $\Lcal_\L=\rho\al$, there exists a single loop with length approximately $\rhoe\al$. In \cite[Theorem 2]{armendariz2011conditional} it is shown that the remaining values are asymptotically iid (in the total variation norm).

It remains to show that the long loop is asymptotically negligible. However, using Lemma \ref{LemAd}, the probability that the single long loop intersects a compact set $K$ is given by $p(x)=o(1)$ and is hence negligible in the evaluation of test functions. This concludes the proof.\qed
\end{proofsect}
\subsection{Proof of Theorem \ref{THM2}, the mean field case}\label{subsec_proof_thm2}
Given the previous results, the proof of the mean field case is not too difficult: the exponential decay of the Hamiltonian forces the total loop length to stay very close to $\rho\al$. We can then apply almost the same proofs as in the free case.


Recall that for $\Delta\subset\Z^d$ bounded,  $\Lcal_\Delta=\abs{\Delta}\Lcalo_\Delta=\sum_{x\in\Delta}\Lcal_x$. We set $\Theta(\eta)=\Theta_\Delta(\eta)=\rho\abs{\Delta}-\Lcal_\Delta(\eta)$, the difference between $\rho\abs{\L}$ and the total particle number. 

Define a new measure $\Pfrak_\Delta$, which \textit{approximates} our limiting distribution, by its Radon--Nikodým derivative
\begin{equation}\label{ApproxPMFDEF}
    \d\Pfrak_\Delta=\frac{1}{\Zfrak_\Delta}\ex^{-a\rho\left[\ell(\omega)-\Theta(\eta)\right]}\1\{\ell(\omega)\ge\Theta(\eta)\}\d\P_\Delta(\eta)\otimes \d M_\Delta(\omega)\, ,
\end{equation}
where
\begin{equation}
    \Zfrak_\Delta=\int \ex^{-a\rho\left[\ell(\omega)-\Theta(\eta)\right]}\1\{\ell(\omega)\ge\Theta(\eta)\}\d\P_\Delta(\eta)\otimes \d M_\Delta(\omega)\, .
\end{equation}
The reference measure $\d\P_\Delta(\eta)\otimes \d M_\Delta(\omega)$ samples a configuration of loops $\eta$ and a single loop $\omega$ independently. $\Pfrak_\Delta$ then only allows configurations which satisfy $\ell(\omega)+\Lcal(\eta)\ge \rho\abs{\L}$ and exponentially penalises configurations by how much $\ell(\omega)+\Lcal(\eta)$ exceeds $\rho\abs{\L}$.

Recall that $\Hpmf_\Delta(\eta)=a\abs{\Delta}\overline{\Lcal}_\Delta(\eta)^2/2$. It admits the following expansion
\begin{equation}\label{ExpansionHam}
    \Hpmf_\Delta(\eta+\delta_\omega)=a\frac{\rho^2\abs{\L}}{2}+a\rho\left[\ell(\omega)-\Theta(\eta)\right]+a\frac{\left[\ell(\omega)-\Theta(\eta)\right]^2}{2\abs{\L}}\, .
\end{equation}
This expansion will help us to absorb the influence of the Hamiltonian into the measure $\d\Pfrak_\Delta$.
\begin{lemma}\label{Equality of partition fuction}
We have that as $N\to\infty$
\begin{equation}
     \E_\L\left[\ex^{- \Hpmf_\L}\1\{\Lcalo_\L\ge \rho\}\right]\sim \Zfrak_\L\ex^{-a\rho^2\abs{\L}/2}\, .
\end{equation}
\end{lemma}
\begin{proofsect}{\textbf{Proof of Lemma \ref{Equality of partition fuction}}}
The idea of the proof is to combine the explicit distribution of $\Lcal_\L$ with the fact that $p_{\rhoe\al}(0)\sim p_{\rhoe\al+j}(0)$ as long as $\abs{j}\ll \al$. This implies that the partition function is not influenced by the CLT fluctuations of $\Lcal_\L$.

Recall that by Corollary \ref{CorAvLoclExceed}, for any $\delta>0$ fixed and $r>\abs{\L}(1+\delta)\rho_\mathrm{c}$ we have 
\begin{equation}\label{Eq_just}
    \P_\L\left(\Lcal_\L=r\right)\sim \abs{\L}\frac{p_{(r-\rho_\mathrm{c}\abs{\L})\beta}(0)}{r-\rho_\mathrm{c}\abs{\L}}\, .
\end{equation}
We have, as $N\to\infty$,
\begin{equation}\label{equationtemp123}
\begin{split}
    \E_\L\left[\ex^{- \Hpmf_\L}\1\{\Lcalo_\L\ge \rho\}\right]&=\sum_{j\ge\rho\abs{\L} }\ex^{-aj^2/(2\abs{\L})} \P_\L\left(\Lcal_\L=j\right)\\
    &\sim \sum_{j\ge\rho\abs{\L} }\frac{\abs{\L}p_{\beta j-\rho_\mathrm{c}\abs{\L}}(0)\ex^{-aj^2/(2\abs{\L})}}{j-\rho_\mathrm{c}\abs{\L}}\\
    &\sim \ex^{-a\rho^2\abs{\L}/2}\sum_{j\ge 0}\rk{\frac{d}{2\pi }}^{d/2}\frac{\abs{\L}\ex^{-a\rho j}}{\beta^{d/2}(j+\rho_\mathrm{e}\abs{\L})^{1+d/2}}\\
    &\sim \ex^{-a\rho^2\abs{\L}/2}\sum_{j\ge {\rho_\mathrm{e}\abs{\L}}}\rk{\frac{d}{2\pi }}^{d/2}\frac{\abs{\L}\ex^{-a\rho (j-{\rho_\mathrm{e}\abs{\L}})}}{\beta^{d/2}j^{1+d/2}}\, .
\end{split}
\end{equation}
Indeed, while the $(1+o(1))$ term in Equation \eqref{Eq_just} may not be uniform in $r$, we can still use the exponential decay:
\begin{equation}
    \P_\L\left(\Lcal_\L=\rho\abs{\L}\right)\sum_{j\ge 0}\ex^{-a\rho j}\frac{\P_\L\left(\Lcal_\L=\rho\abs{\L}+j\right)}{\P_\L\left(\Lcal_\L=\rho\abs{\L}\right)}\le C \P_\L\left(\Lcal_\L=\rho\abs{\L}\right)\sum_{j\ge 0}\ex^{-a\rho j}\left(1+\frac{j}{\rho_\mathrm{e}\abs{\L}}\right)^{-1-d/2}\, ,
\end{equation}
where the bound $\P_\L\left(\Lcal_\L=\rho\abs{\L}+j\right)\le C\abs{\L}(j+\rho_\mathrm{e}\abs{\L})^{-1-d/2}$ for all $j\ge 0$ is from Proposition \ref{PropositionLDP}. Thus, only finitely many terms in the expansion are relevant, and our $(1+o(1))$ argument in Equation \eqref{equationtemp123} is justified.

Note that we can rewrite Equation \eqref{equationtemp123} as
\begin{equation}
     \E_\L\left[\ex^{- \Hpmf_\L}\1\{\Lcalo_\L\ge \rho\}\right]\sim \ex^{-a\rho^2\abs{\L}/2}M_\L\left[\ex^{-a\rho (\ell(\omega)-{\rho_\mathrm{e}\abs{\L}})}\1\{\ell(\omega)\ge {\rho_\mathrm{e}\abs{\L}}\}\right]\, ,
\end{equation}
where the second use of $\sim$ follows from Equation \eqref{equationtemp123}.

We now show that
\begin{equation}\label{Equation19423}
    \Zfrak_\L\sim M_\L\left[\ex^{-a\rho (\ell(\omega)-{\rho_\mathrm{e}\abs{\L}})}\1\{\ell(\omega)\ge {\rho_\mathrm{e}\abs{\L}}\}\right]\sim\abs{\L}^{-d/2}\sum_{j\ge {\rho_\mathrm{e}\abs{\L}}}\rk{\frac{d}{2\pi }}^{d/2}\frac{\rho_\mathrm{e}^{-1-d/2}\ex^{-a\rho (j-{\rho_\mathrm{e}\abs{\L}})}}{\beta^{d/2}(1+j/\rho_\mathrm{e}\abs{\L})^{1+d/2}}\, .
\end{equation}
Note that the sum on the right-hand side of the above equation is bounded from above and below by a constant.

We first expand to isolate the main contribution
\begin{equation}
    \Zfrak_\L=\Zfrak_{\L,m}+\Zfrak_{\L,-}+\Zfrak_{\L,+}\, .
\end{equation}
Here, set
\begin{equation}
    \Zfrak_{\L,m}=\E_\L\otimes M_\L\left[\ex^{-\rho\left[\ell(\omega)-\Theta(\eta)\right]}\1\{\ell(\omega)\ge\Theta(\eta)\}\1\{\abs{\Lcal_\L-\rho_\mathrm{c}\abs{\L}}\le \abs{\L}^{5/6}\}\right]\, .
\end{equation}
However, on the event $\abs{\Lcal_\L(\eta)-\rho_\mathrm{c}\abs{\L}}\le \abs{\L}^{5/6}$, we have that $\abs{\Theta(\eta)-\rhoe\al}\le \al^{5/6}$ and hence, as $N\to\infty$
\begin{equation}
    M_\L\left[\ex^{-\rho\left[\ell(\omega)-\Theta(\eta)\right]}\1\{\ell(\omega)\ge\Theta(\eta)\}\right]\sim M_\L\left[\ex^{-a\rho (\ell(\omega)-{\rho_\mathrm{e}\abs{\L}})}\1\{\ell(\omega)\ge {\rho_\mathrm{e}\abs{\L}}\}\right]\, ,
\end{equation}
since uniformly for $-\al^{5/6}\le j\le\al^{5/6}$ we have (compare Equation \eqref{equationtemp123}), as $N\to\infty$
\begin{equation}
    M_\L\left[\ex^{-a\rho (\ell(\omega)-{\rho_\mathrm{e}\abs{\L}})}\1\{\ell(\omega)\ge {\rho_\mathrm{e}\abs{\L}}\}\right]\sim M_\L\left[\ex^{-a\rho (\ell(\omega)-{\rho_\mathrm{e}\abs{\L}}-j)}\1\{\ell(\omega)\ge {\rho_\mathrm{e}\abs{\L}}+j\}\right]\, .
\end{equation}
Furthermore, note that $\P_\L\left(\abs{\Lcal-\rho_\mathrm{c}\abs{\L}}\le \abs{\L}^{5/6}\right)=1+o(1)$. This gives that $\Zfrak_{\L,m}$ is asymptotically equal to the right-hand side of Equation \eqref{Equation19423}.

Next, we show that the case where $\abs{\Lcal-\rho_\mathrm{c}\abs{\L}}> \abs{\L}^{5/6}$ is negligible. We set
\begin{equation}
    \Zfrak_{\L,-}=\E_\L\otimes M_\L\left[\ex^{-\rho\left[\ell(\omega)-\Theta(\eta)\right]}\1\{\ell(\omega)\ge\Theta(\eta)\}\1\{\Lcal_\L-\rho_\mathrm{c}\abs{\L}\le -\abs{\L}^{5/6}\}\right]\, .
\end{equation}
However, as $\Zfrak_{\L,0}$ is of polynomial order and $\P_\L\left(\Lcal_\L-\rho_\mathrm{c}\abs{\L}\le -\abs{\L}^{5/6}\right)$ decays at a stretch-exponential speed (recall Proposition \ref{PropositionLowerBound}), we can neglect it.

For $\Zfrak_{\L,+}=\Zfrak_\L-(\Zfrak_{\L,0}+\Zfrak_{\L,-})$, a similar argument implies that it is negligible. This concludes the proof.\qed
\end{proofsect}

The next proposition is a summary of some technical results, the analogy of the technical lemmas for the free case.
\begin{proposition}\label{PropTechnicalMeanField}
For $K\subset\Z^d$ bounded, we have as $N\to\infty$
\begin{align}
        &\P_N^{\pmf,\rho}\left(\exists \omega_o\colon \omega(0)\in\L_N,\, \omega_o\cap K\neq\emptyset,\,\forall\omega\colon \ell(\omega)\notin\Re\right)=o(1)\, ,\\
         &\P_N^{\pmf,\rho}\left(\exists x\in C_N\setminus\gk{0}\colon \textnormal{two or more loops started at }xN+\L \textnormal{ intersect }K\right)=o(1)\, ,\\
        &\P_N^{\pmf,\rho}\left(\exists \omega\colon \omega(0)\notin C_N\,\,\mathrm{ but }\,\,\omega\cap K\neq\emptyset\right)=o(1)\, ,\\
        &\P_N^{\pmf,\rho}\left(\exists\omega_o\colon \omega_o(0)\notin\L,\,\omega_o\cap K\neq\emptyset,\,\ell(\omega_o)\notin\Re \right)=o\left(1\right)\, ,\\
        &\P_N^{\pmf,\rho}\left(\exists\omega_o\colon \omega_o(0)\in\L,\,\omega_o\cap K\neq\emptyset,\,\ell(\omega_o)\in\Re \right)=o\left(1\right)\, .
\end{align}
\end{proposition}

\begin{proofsect}{\textbf{Proof of Proposition \ref{PropTechnicalMeanField}}}
The proof of the above statements is the same as in the free case, because the Hamiltonian $\Hpfm_\Delta$ is determined by the value of $\Lcal_\Delta$. We illustrate this with the first statement.

The proof of the first statement is similar to the one of Lemma \ref{LemmaAtLeastOneLongLoop}. Note that using Lemma \ref{Equality of partition fuction} in the last step, for $\Delta=xN+\L$
\begin{equation}
\begin{split}
     \sum_{x\in C_N}&\Ppmf_{xN+\L}\left(\exists \omega_o\colon \omega_o\cap K\neq\emptyset,\,\forall\omega\colon \ell(\omega)\notin\Re\,\Big|\Lcalo_\Delta\ge \rho\right)\\
     &=\sum_{x\in C_N}\frac{\E_{xN+\L}\left[\ex^{-\Hpfm_\Delta},\,\exists \omega_o\colon \omega_o\cap K\neq\emptyset,\,\forall\omega\colon \ell(\omega)\notin\Re\, ,\Lcalo_\Delta\ge \rho\right]}{\E_\L\left[\ex^{- \Hpmf_\L}\1\{\Lcalo_\Delta\ge \rho\}\right]}\\
     &\sim\sum_{x\in C_N}\frac{\E_{xN+\L}\left[\ex^{-\Hpfm_\Delta},\,\exists \omega_o\colon \omega_o\cap K\neq\emptyset,\,\forall\omega\colon \ell(\omega)\notin\Re\, ,\Lcalo_\Delta\ge \rho\right]}{\Zfrak_\L\ex^{-a\rho^2\abs{\L}/2}}\, .
     \end{split}
\end{equation}
We only need to bound the denominator: we expand
\begin{equation}\label{EquationSum4114}
\begin{split}
    \E_{xN+\L}&\left[\ex^{-\Hpfm_\Delta},\,\exists \omega_o\colon \omega_o\cap K\neq\emptyset,\,\forall\omega\colon \ell(\omega)\notin\Re,\,\Lcalo_\Delta\ge \rho_\mathrm{c}+\rho_\mathrm{e}\right]\\
    &=\sum_{k\ge \rho\abs{\L}}\ex^{-\Hpmf( k)}\P_{xN+\L}\left(\exists \omega_o\colon \omega_o\cap K\neq\emptyset,\,\forall\omega\colon \ell(\omega)\notin\Re,\,\Lcal_\Delta=k\right)\, ,
\end{split}
\end{equation}
where we wrote $\Hpmf( k)$ for $\Hpmf_\Delta(\eta)$ evaluated at $\Lcal_\Delta(\eta)=k$.

Lemma \ref{LemmaAtLeastOneLongLoop} shows that there is some $\delta>0$ such that for all $k$ with $\rho\al\le k\le \rho\al+\al^{\delta}$ we have that, as $N\to\infty$
\begin{equation}
    \sum_{x\in C_N}\P_{xN+\L}\left(\exists \omega_o\colon \omega_o\cap K\neq\emptyset,\,\forall\omega\colon \ell(\omega)\notin\Re,\,\Lcal_\Delta=k\right)=o\rk{\al^{-d/2}}=o\left(\Zfrak_\L \right)\, .
\end{equation}
However, for $k>\rho\al+\al^{\delta}$, we have that $\ex^{-\Hpmf( k)}<\ex^{-a\rho^2\al/2}\ex^{-a\al^{\delta}}$ and they hence do not contribute to the sum in Equation \eqref{EquationSum4114}.

This concludes the proof of the first statement. The remaining statements follow in the same way: we can transfer the proof from the free case by conditioning on the value of $k$ and noting that $k$ larger than $\rho\al$ are suppressed exponentially fast.\qed
\end{proofsect}
Note that by the last statement of Proposition \ref{PropositionPoisson} and Equation \eqref{ExpansionHam}
\begin{equation}
    \E_{xN+\L}\left[\ex^{-H^\pmf_\Delta}F,\exists! \omega \colon\ell(\omega)\in\Re,\,\Lcal_\Delta\ge \rho\abs{\L}\right]\sim \EPfrak_{xN+\L}\left[\ex^{-H^\pmf_\Delta}F\ex^{a\rho\left[\ell(\omega_o)-\Theta(\eta)\right]}\right]
\end{equation}
Using Proposition \ref{PropTechnicalMeanField} in the same way as we used it in the proof of Proposition \ref{LemStTry}, we obtain
\begin{equation}
    \E_N^{\pmf,\rho}[F]\sim\bigotimes_{x\in C_N}\EPfrak_{xN+\L}\left[F\ex^{-\Hpfm_\Delta+a\rho^2\abs{\L}/2+\rho\left[\ell(\omega_o)-\Theta(\eta)\right]}\right]\, .
\end{equation}
Next, we want to remove the influence of the Hamiltonian. For this, we add and subtract $1$ from $\ex^{-\Hpfm_\Delta+a\rho^2\abs{\L}/2+\rho\left[\ell(\omega_o)-\Theta(\eta)\right]}$. We first estimate the error term:
\begin{lemma}\label{RemoveHamil}
We have that for $K\subset \Z^d$ compact, as $N\to\infty$
\begin{equation}\label{Equation194231}
    \sum_{x\in C_N}\EPfrak_{xN+\L}\left[\Big|\ex^{-\Hpfm_\Delta+a\rho^2\abs{\L}/2+\rho\left[\ell(\omega_o)-\Theta(\eta)\right]}-1\Big|,\,\omega\cap K\neq \emptyset\right]=o(1)\, .
\end{equation}
\end{lemma}
\begin{proofsect}{\textbf{Proof of Lemma \ref{RemoveHamil}}}
The proof is similar to the second part of the proof of Lemma \ref{LemAd}, hence we focus on the steps different and refer back for additional details.

Assume without loss of generality that $K$ is the origin $K=\gk{0}$. We multiply Equation \eqref{Equation194231} by $\Zfrak_\L$ and expand the probability (similar to the previous lemmas)
\begin{multline}
 \sum_{k\ge 0}\sum_{j\ge \rho\al-k}
    \sum_{x\in C_N}\P_{xN+\L}\otimes M_{xN+\L}\left[\Big|\ex^{-\Hpfm_\Delta+a\rho^2\abs{\L}/2+\rho\left[\ell(\omega_o)-\Theta(\eta)\right]}-1\Big|,\omega\cap K\neq \emptyset,\, \Lcal_\L=k,\ell=j\right]\\
    \le \sum_{k\ge 0}\sum_{j\ge \rho\al-k}\P_\L\left(\Lcal_\L=k\right)M_{\L_N}\left[\ell=j,0\in\omega\right]\Babs{\ex^{-\Hpfm_\Delta(k+j)+a\rho^2\abs{\L}/2}-\ex^{-a\rho\left[j-\rho\al+k\right]}}\, ,
\end{multline}
where we abbreviate $a(k+j)^2/2\al$ by $\Hpfm_\Delta(k+j)$.

Recall that $M_{\L_N}\left[\ell=j,0\in\omega\right]=\Ocal\rk{j^{-d/2}}$, see Equation \eqref{EquationEstimateByRange}. Note that using Equation \eqref{ExpansionHam}
\begin{equation}
    -\Hpfm(k+j)+a\rho^2\abs{\L}/2=-a\rho[j-k]-a\frac{[j-k]^2}{2\abs{\L}}\, .
\end{equation}
By re-indexing the sum, we rewrite this as
\begin{equation}
    \Ocal\left(1\right) \sum_{k\ge 0}\sum_{j\ge 0}\P_\L\left(\Lcal_\L=k\right)\left(j+\rho\al-k\right)^{-d/2}\Babs{\ex^{-j^2/2\abs{\L}}-1}\ex^{-a\rho j}\, .
\end{equation}
As $k$ concentrates around $\abs{\L}\rho_\mathrm{c}$ and $j$ around $\rhoe\al$, we get that $\ex^{-j^2/2\abs{\L}}-1=o(1)$ for contributing $j$. Hence, we can estimate the above as $o\left(\abs{\L}^{-d/2}\right)$. As $\Zfrak_\L=\Ocal\left(\abs{\L}^{-d/2}\right)$, this proves the claim.\qed
\end{proofsect}
The next lemma is the analogue of Proposition \ref{LemStTry}.
\begin{lemma}\label{LemmaConvRandInt}
We define for $\Delta=xN+\L$, the measure $\Mfrak_\Delta$ as the projection of $\Pfrak_\Delta$ onto its second component, the long loop $\omega$.

We furthermore set for $K\subset \Z^d$
\begin{equation}
    \d{\Mfrak_K}=\sum_{x\in C_N}\1\{\omega\cap K\neq\emptyset\}\d\Mfrak_\Delta\, .
\end{equation}
Define ${\Mfrak_K}^{*}={\Mfrak_K}\circ\Pi^{-1}$, analogous to Proposition \ref{LemStTry}. We then have that for $E$ an events as in Definition \ref{DefApprx}
\begin{equation}
    {\Mfrak_K}^{*}[E]=\rho_\mathrm{e}\nu[E]\left(1+o(1)\right)\, ,
\end{equation}
as $N\to\infty$
\end{lemma}
\begin{proofsect}{\textbf{Proof of Lemma \ref{LemmaConvRandInt}}}
Following the steps from the proof of Proposition \ref{LemStTry}, we can rewrite for $z$ a point in the boundary of $K$, as $N\to\infty$
\begin{equation}\label{Equation52623}
\begin{split}
    {\Mfrak_K}^{*}[E]&\sim \kappa_d q(E) {\Mfrak_K}[\ell(\omega)\1\{\omega(0)=z\}]\\
    &\sim \frac{\kappa_d q(E)}{\Zfrak_\L}\sum_{j=-\al^{5/6}}^{\al^{5/6}}\P_\L\rk{\Lcal_\L=\rhoc\al+j}\sum_{k\ge \rhoe\al-j}\rk{\frac{d}{2\pi }}^{d/2}\frac{\ex^{-a\rho(\beta k-\rhoe\al+j)}}{(\beta k)^{d/2}}\, .
    \end{split}
\end{equation}
where $q(E)$ is the same as in Equation \eqref{EqQE} and we used Lemma \ref{Equality of partition fuction} to neglect the event $\gk{\abs{\Lcal_\L-\rhoc\al}>\al^{5/6}}$. We rewrite for $j$ in the range as above, as $N\to\infty$
\begin{equation}
\begin{split}
   \sum_{k\ge \rhoe\al-j}\frac{\ex^{-a\rho(\beta k-\rhoe\al+j)}}{(\beta k)^{d/2}}&=\frac{1}{\rk{\rhoe\al(1+o(1))}^{d/2}}\sum_{k\ge 0}\frac{\ex^{-a\rho k  }}{(1+\beta k/\rhoe\al)^{d/2}}\\
    &=\frac{\left(1+o(1)\right)}{(\rho_\mathrm{e}\abs{\L})^{d/2}\left(1-{\ex^{-a\rho}}\right)}\, .
\end{split}
\end{equation}

Similarly, we find that
\begin{equation}
   \Zfrak_\L\sim\sum_{j\ge {\rho_\mathrm{e}\abs{\L}}}\rk{\frac{d}{2\pi }}^{d/2}\frac{\abs{\L}\ex^{-a\rho (j-{\rho_\mathrm{e}\abs{\L}})}}{\beta^{d/2}j^{1+d/2}}\sim\rk{\frac{d}{2\pi }}^{d/2} \frac{\abs{\L}}{(\rho_\mathrm{e}\abs{\L})^{d/2+1}\left(1-{\ex^{-a\rho}}\right)}\, .
\end{equation}
Therefore, we can conclude that
\begin{equation}
    {\Mfrak_K}^{*}[E]=\rho_\mathrm{e}\nu[E]\left(1+o(1)\right)\, .
\end{equation}
This concludes the proof.\qed
\end{proofsect}
The proof of Theorem \ref{THM2} now goes as follows: the distribution of $\Ppmf_N$ conditioned on exceeding a particle density of $\rho$ can be obtained in two steps. First, sample $\eta=\sum_{\omega}\delta_\omega$ according to $\P_N$. Then, sample a PPP process with intensity measure $\Mfrak_K$. However, this has the same distribution as a PPP with intensity measure $\nu$, according to Lemma \ref{LemmaConvRandInt}: take $F$ a suitable test function (compactly supported). Then, using Proposition \ref{PropTechnicalMeanField} in the first step and
\begin{equation}
    \begin{split}
         \E_N^\pmf[F]&\sim\bigotimes_{x\in C_N}\EPfrak_{xN+\L}\left[F\ex^{-\Hpfm+a\rho^2\abs{\L}/2+\rho\left[\ell(\omega_o)-\Theta(\eta)\right]}\right]\\
         &\sim \bigotimes_{x\in C_N}\EPfrak_{xN+\L}\left[F\right]+\Ocal\left(\bigotimes_{x\in C_N}\EPfrak_{xN+\L}\left[\Big|1-\ex^{-\Hpfm+a\rho^2\abs{\L}/2+\rho\left[\ell(\omega_o)-\Theta(\eta)\right]}\Big|\1\{\eta\cap\{0\}\neq \emptyset\}\right]\right)\\
         &\sim \E_{\Z^d,\beta,0}\otimes\textsf{E}_{{\Mfrak_K}^{*}}[F]\sim \E_{\Z^d,\beta,0}\otimes \E_{\rho-\rho_\mathrm{c}}^\iota\left[F\right]\, ,
    \end{split}
\end{equation}
where we used the notation $\textsf{E}_{{\Mfrak_K}^{*}}$ for the PPP with intensity measure ${\Mfrak_K}^{*}$, just like in Proposition \ref{LabelLemmaMain}. This concludes the proof.\qed
\section*{Appendix}
\subsection*{Proof of Theorem \ref{ThMWholeDescription}}
To make this paper self-contained, we give the proof of the case $\rho<\rhoc$ here. It can also be found in \cite{dickson2021formation}.

Fix $\rho<\rhoc$. Since the function $\boldsymbol{\rho}$ (from Equation \eqref{Definitionbrho}) is bijective from $(-\infty,0]\to (0,\rhoc]$, we can choose $\mu$ such that $\brho(\mu)=\rho$. We now use Equation \eqref{EquationChangeOfMeasure} to write for any test-function and any $\tilde{\mu}$
\begin{equation}
    \E_{\L,\beta,\tilde{\mu}}[F|\Lcal_\L=\rho\al]=\E_{\L,\beta,{\mu}}[F|\Lcal_\L=\rho\al]\, .
\end{equation}
By the Campbell formula (Proposition \ref{PropositionPoisson}), we have that $\E_{\L,\beta,{\mu}}[\Lcal_\L]=\brho(\mu)\al=\rho\al$. By Gnedenko’s Local Limit Theorem (see \cite[Theorem 8.4.1]{bingham1989regular}), we have that for any $M>0$ fixed 
\begin{equation}
    \sup_{j\colon \abs{j-\rho\al}\le Ma_\L}\abs{a_\L\P_{\L,\beta,{\mu}}\rk{\Lcal_\L=j}-\phi\rk{\frac{j-\rho\al}{\a_\L}}}=o(1)\, ,
\end{equation}
where $\phi$ is the density of the limiting stable law, compare Corollary \ref{CorAvLoclExceed}. Hence, for all $\delta>0$ such that $\al^\delta=o\rk{a_\L}$, for all $-\al^\delta\le j\le \al^\delta$
\begin{equation}\label{EquationSep}
    \P_{\L,\beta,{\mu}}\rk{\Lcal_\L=\rho\al}\sim \P_{\L,\beta,{\mu}}\rk{\Lcal_\L=\rho\al+j}\, .
\end{equation}
We now can proceed similar to the supercritical case: fix a region $\Delta_M=[-M,M)^d\cap\Z^d$ for $M=N^{\delta}$. Fix a compact set $K\subset\Z^d$. The loops intersecting $K$ start in $\Delta_M$ with high probability, see Lemma \ref{LemAd}. With high probability, the value of $\Lcal_{\Delta_M}$ is bounded by $2\rho\abs{\Delta_M}$. Equation \eqref{EquationSep} gives that the distribution of $\Lcal_\L$ is not affected by $\Lcal_{\Delta_M}$ and hence we can write for test functions $F$
\begin{equation}
    \E_{\L,\beta,{\mu}}[F|\Lcal_\L=\rho\al]\le  \E_{\Delta_M,\beta,{\mu}}[F,\Lcal_{\Delta_M}\le 2\rho\abs{\Delta_M}]\sup_{j\colon \abs{j-\rho\al}\le Ma_\L}\frac{\P_{\L,\beta,{\mu}}\rk{\Lcal_\L=\rho\al+j}}{\P_{\L,\beta,{\mu}}\rk{\Lcal_\L=\rho\al}}\, ,
\end{equation}
and similar for a lower bound. Equation \eqref{EquationSep} then shows that only $\E_{\Delta_M,\beta,{\mu}}[F,\Lcal_{\Delta_M}\le 2\rho\abs{\Delta_M}]$ survives in the limit. Since $\E_{\Delta_M,\beta,{\mu}}[\Lcal_{\Delta_M}]\rho\abs{\Delta_M}]$, we have that
\begin{equation}
    \E_{\Delta_M,\beta,{\mu}}[F,\Lcal_{\Delta_M}\le 2\rho\abs{\Delta_M}]\sim \E_{\Z^d,\beta,{\mu}}[F]\, ,
\end{equation}
as $N\to\infty$. This concludes the proof.\qed

\subsection{Table with frequently used notation}
\mbox{}
  
\begin{center}
{  \renewcommand{\arraystretch}{1.4}
 \captionof{table}{List of frequently used notation}
\begin{tabular}{ |p{1.4cm}||p{6cm}|p{4.2cm}|p{2.85cm}|  }
 \hline
 Symbol & Definition & Explanation &Class\\
 \hline \hline
 $\rho$   &  Positive real  &Density&   Model parameter\\
 \hline
  $\mu$   &   Non-positive  &Chemical potential&   Model parameter\\
 \hline
 $\beta$   &  Positive real   &Inverse temperature&   Model parameter\\
 \hline
 $\rho_\mathrm{c}$   &  $\sum_{j\ge 1}p_{\beta j}(0)$   &Critical density&   Parameter\\
 \hline
 $\rho_\mathrm{e}$   &  $\max\{\rho-\rho_\mathrm{c},0\}$   &Excess density&   Parameter\\
 \hline
 $\L$   &  $[-N/2,N/2)^d\cap\Z^d$   & Box &   Set\\
  \hline
   $\L_N$   &  $r_N[-N^{d/2}/2,N^{d/2}/2)^d\cap\Z^d$  & Box, Section \ref{subsectionDomain}  &   Set\\
  \hline
     $C_N$   &  $r_N[-N^{d/2-1}/2,N^{d/2-1}/2)^d\cap\Z^d$   & Box &   Set\\
  \hline
     $\Re$   &  $[\rhoe\al-\al^{5/6},\rhoe\al+\al^{5/6}]\cap \Z$   & Interval &   Set\\
  \hline
 $M_{\L,\beta,\mu}$   &  See Equation \eqref{FirstEquationOfM}   &Bosonic Loop measure &   Measure\\
  \hline
 $ \overline{M}_{\Delta}$   &  See Equation \eqref{Mbardef}   & Long loop measure &   Measure\\
  \hline
 $ Z_{\Delta}$   &  See Equation \eqref{EquationPartFuncDef}   & Partition function &   Measure\\
  \hline
 $\P_{\L,\beta,\mu}$   &  PPP with intensity measure  $M_{\L,\beta,\mu}$  & Loop soup &   Measure\\
   \hline 
 $\P_{N}^\rho$   & See Equation \eqref{EquationDefinitionrhon} & Conditional process  &   Measure\\
   \hline  
$\P^\pmf_{\Delta}$   & See Equation \eqref{EquationDefineThe MFmodel} & Mean field law  &   Measure\\
   \hline  
$\Pfrak_\Delta$   & See Equation \eqref{ApproxPMFDEF} & Mean field approx.  &   Measure\\
   \hline  
 $\P_{\rhoe}^\iota$   &  See Section \ref{sectionInterlacements}  & Interlacements &   Measure\\
    \hline  
 $\P_{x,x}^{\beta j}$   &  See Equation \ref{EquationOfBridgeMeasure}  & Random walk bridge &   Measure\\
     \hline  
 $p_{\beta j}(x,y)$   &  See Section \ref{subsectionDomain}  & transition kernel &   Kernel\\
      \hline  
 $\Lcal_\Delta$   &  See Equation \ref{DefinitionLocalTime}  & Particle number &   Random variable\\
       \hline  
 $\Lcalo_\Delta$   &  $\abs{\Delta}^{-1}\Lcal_\Delta$  & Particle density &   Random variable\\
 \hline
\end{tabular}\label{table1} }
\end{center}

\subsection*{Acknowledgements} The author would like to thank the anonymous referees for their many helpful suggestions. The author would also like to thank Julius Damarackas, Chokri Manai and Silke Rolles for their help regarding typos and presentation.
\bibliography{thoughts}{}
\bibliographystyle{alpha}
\end{document}